\documentclass[reqno,11pt]{article}

\usepackage[utf8]{inputenc}
\usepackage[T1]{fontenc}
\usepackage{amsmath}
\usepackage{amssymb}
\usepackage{amsthm}
\usepackage{mathrsfs}
\usepackage{xcolor}
\usepackage{graphicx}
\usepackage{tikz}
\usetikzlibrary{patterns, arrows.meta}
\usepackage{comment}
\usepackage[numbers]{natbib}
\usepackage{hyperref}
\usepackage{geometry}
\newtheorem{theorem}{Theorem}[section]
\newtheorem{lemma}[theorem]{Lemma}

\theoremstyle{definition}
\newtheorem{definition}[theorem]{Definition}

\theoremstyle{plain}
\newtheorem{corollary}[theorem]{Corollary}
\newtheorem{proposition}[theorem]{Proposition}

\theoremstyle{remark}
\newtheorem{remark}[theorem]{Remark}

\numberwithin{equation}{section}

\title{Determination of a Bi-periodic Refractive Index for the Helmholtz Equation}
\author{ 
	Han Chen\thanks{School of Mathematical Sciences, Nankai University, Tianjin 300071, China. Email: \texttt{hanc@mail.nankai.edu.cn}}\quad and
	Guanghui Hu\thanks{(Corresponding author) School of Mathematical Sciences and LPMC, Nankai University, Tianjin 300071, China. Email: \texttt{ghhu@nankai.edu.cn}} 
}
\date{}

\begin{document}
	\maketitle
	
	\begin{abstract}
		We investigate an inverse problem for the Helmholtz equation in a three-dimensional half-space with a biperiodic medium above a
		sound-soft plane. We establish the uniqueness of a real-valued, bounded, biperiodic refractive index that depends only on the horizontal variables using spectral techniques. 
		 An orthogonality identity is first derived from near-field measurements. By constructing suitable separated solutions, this
		 identity is reduced to relations involving eigenfunctions of
		 two-dimensional quasiperiodic Schr\"odinger operators. To recover
		 Fourier information from these relations, we project horizontal
		 plane waves onto suitable high-energy spectral subspaces, obtaining
		 finite combinations of eigenfunctions that approximate the plane
		 waves with a quantified error.  Combining these approximations
		with the orthogonality identity, we show that all Fourier coefficients of the difference of two refractive indices vanish, thereby establishing the desired uniqueness.
		\vspace*{0.5cm}
		
		\vspace*{0.5cm}
		
	\end{abstract}
	
	\section{Introduction}

The reconstruction of the refractive index of a periodic or bi-periodic medium from measurements of the scattered wave field is a central inverse problem of wave physics, and it underlies a broad class of applications in modern science and technology \cite{BaoDobsonCox1995,BL2025,Petit1980}. In micro-optics and nanophotonics, diffraction gratings and photonic crystals are engineered periodic structures whose optical response, including diffraction efficiency, dispersion and photonic band gaps, is governed entirely by the spatial distribution of the refractive index.  Thus, recovering this distribution from measured data is the key step in their non-destructive characterization, quality control and inverse design \cite{BL2025,Petit1980}. In semiconductor manufacturing, so-called scatterometry (optical critical-dimension metrology) reconstructs the geometry and the optical constants of periodic photoresist gratings from their near-field or reflected response. In optical holography the recorded interference pattern is stored precisely as a periodic modulation of the refractive index, and reading the hologram amounts to inverting this modulation. Related reconstruction problems arise in the design and testing of anti-reflection coatings and photonic-crystal fibers and sensors, in acoustic scattering by periodic structures and inhomogeneous layers, where the index describes the sound-speed profile of absorbing liners and of rough or layered interfaces \cite{CMT,Zhang B2018}, and in radar and remote sensing of layered periodic media. In all of these settings, only a limited amount of information is accessible, typically  the field on a measurement surface above the structure. This naturally raises the basic mathematical question of the extent to which the unknown index is uniquely determined by such data.  The present paper is devoted to this uniqueness question for a bi-periodic medium in three dimensions.

The direct problem of time-harmonic scattering by periodic structures and layered media has been investigated extensively over the last three decades \cite{BaoDobsonCox1995,BBS1994,Kirsch A1995,Petit1980}. In the electromagnetic case, variational frameworks for diffraction by general bi-periodic structures in $\mathbb{R}^{3}$ establish the existence and uniqueness of quasiperiodic solutions and provide the basis for the rigorous analysis and numerical approximation of such problems \cite{BaoDobsonCox1995,BL2025}. For acoustic scattering by periodic inhomogeneous media, and for inhomogeneous layers on rigid plates, related variational  formulations were investigated in \cite{ BBS1994, CMT,  ElschnerSchmidt1998, HuLiuQuZhang2015,Strycharz B 1998, Zhang B2018}. More recently, Kirsch and Lechleiter \cite{Kirsch A2018} analyzed scattering by penetrable periodic layers in the upper half-space, proving a limiting absorption principle and formulating a Rayleigh-type radiation condition that remains valid in the presence of surface waves and guided modes; related results for locally perturbed periodic structures and for bi-periodic inhomogeneous layers were obtained in \cite{HuKi26,HuKiZh}. We also mention that, since the Rayleigh radiation condition alone does not always guarantee uniqueness, Bonnet-Bendhia and Starling \cite{BBS1994} constructed explicit examples of guided (surface) waves that lead to non-uniqueness for the diffraction problem. These works provide the well-posedness framework on which the present study relies.

Inverse problems for periodic gratings and layers have likewise attracted considerable attention. For electromagnetic diffraction gratings, uniqueness results for determining the profile of perfectly conducting or penetrable periodic surfaces have been obtained by a variety of methods, including Rellich-type identities and orthogonality arguments \cite{Kirsch1994,Kirsch A1995}. Hu, Yang and Zhang \cite{hu2011} considered Maxwell's equations in a bi-periodic inhomogeneous layer on a perfectly conducting plate and proved that a refractive index $q$ which depends only on $x_{1}$ inside the layer is uniquely determined by the tangential electric field on a plane above the layer; their proof combines an orthogonality relation for two refractive indices with spectral information from a quasiperiodic Sturm--Liouville eigenvalue problem in the vertical direction. Gibson, Hu and Zhao \cite{Gibson Hu and Zhao} studied inverse medium scattering from infinite biperiodic structures in three dimensions and established uniqueness from multi-frequency plane-wave data, again exploiting Floquet--Bloch expansions and orthogonality relations. Strycharz \cite{Strycharz B.} established uniqueness in the inverse transmission problem for periodic media.

The purpose of this paper is to recover a real-valued, bounded, bi-periodic refractive index $q=q(x_{1},x_{2})$ in the periodic cell from near-field measurements of the total field on a plane above the structure. In the corresponding inverse problem for a bounded inhomogeneous medium, one may combine an integral identity with complex geometrical optics solutions whose phases can be chosen to probe arbitrary Fourier modes of the unknown coefficient\cite{Sylvester J1987}. In the present setting, however, the admissible particular solutions are required to be $\alpha$-quasiperiodic, which prevents the direct use of such a construction. Our approach therefore follows the strategy pioneered by Kirsch \cite{Kirsch A1995} and later extended by Hu, Yang and Zhang \cite{hu2011} in the electromagnetic case: we first derive an orthogonality relation for the difference $q_{1}-q_{2}$ of two refractive indices that generate identical near-field data on the measurement plane,  and then construct suitable interior
solutions to recover its Fourier coefficients. Since the indices
are independent of the vertical variable inside the layer, the
construction leads to the two-dimensional quasiperiodic
Schr\"odinger operators
\[
 L(\tilde{q}_j):=-\Delta_{\tilde{x}}+\tilde{q}_j(\tilde{x}),\quad \tilde{q}_j(\tilde{x}):=-k^{2}q_j(\tilde{x}),\quad\tilde{x}=(x_1,x_2)\in(0,2\pi)^{2},
\]
endowed with the prescribed $\alpha$-quasiperiodic boundary
conditions. For real-valued bounded potentials, the  Bloch fiber operators $L_{t}(\tilde{q}_j)$ are self-adjoint on $L^{2}(F)$ with compact resolvent, and therefore have purely discrete spectra.

It is instructive to compare our problem with the one treated by Kirsch \cite{Kirsch A1995}. There the structure is periodic in a single direction and the unknown index depends on $x_{1}$ only; consequently, the particular solutions can be obtained by separation of variables and the whole analysis is reduced to a one-dimensional quasiperiodic Sturm--Liouville problem. In that one-dimensional setting, for a fixed quasi-momentum
$\tau\notin\frac12\mathbb{Z}$, every sufficiently large eigenvalue
is simple and well separated from the remaining spectrum, and the
corresponding eigenfunction is dominated by a single plane wave. Hence the eigenpairs admit asymptotic expansions of the form
\[
\lambda_{n}=(n+\tau)^{2}-\frac{k^2}{2\pi}\int_0^{2\pi}q(s)\,ds+o(1),\qquad \psi_{n}(x_{1})=e^{i(n+\tau)x_{1}}+O(1/n),\qquad n\to\infty.
\]
 These expansions are obtained by elementary perturbation theory, or even directly from the fundamental solutions and the characteristic equation. Substituting them into the orthogonality relation and letting $n\to\infty$,  then recovers each Fourier coefficient of $q_{1}-q_{2}$.

The bi-periodic problem studied here is essentially more difficult, because the reduction above is no longer available and the inverse problem is governed by a genuinely two-dimensional periodic Schr\"odinger operator. Moreover, distinct lattice points may
give equal or closely spaced free eigenvalues $|n+\alpha|^2$,
and the potential may couple the corresponding Fourier modes\cite{V2024}. 
Consequently, an approximation of each prescribed plane wave by
an individual eigenfunction cannot be assumed. What the
orthogonality argument requires, however, is a family of spectral
combinations whose products recover prescribed Fourier modes.
This motivates the use of spectral projections. We project each
$\alpha$-quasiperiodic plane wave onto a suitable high-energy
spectral interval of the corresponding operator. A uniformly
bounded residual gives an $L^2(F)$ approximation error of order
$|n+\alpha|^{-1}$. The projection retains entire eigenspaces and
therefore provides the required plane-wave approximation at the
fixed parameter $\alpha$, including parameters at which spectral
multiplicities occur. 

To carry out this programme, we use an $\alpha$-quasiperiodic point
source as the incident field and measure the total field for all
source and receiver points on $\Gamma_a$. We first prove that the
Robin traces of these total fields on the upper boundary of the
layer span a dense subspace of the relevant trace space. The coincidence
of the near-field data then implies the coincidence of  interior
Robin-to-Dirichlet maps, from which the orthogonality relation for
$q_1-q_2$ follows. We next construct separated interior solutions
by combining horizontal eigenfunctions with vertical profiles
satisfying the Dirichlet boundary condition. For every pair of
positive horizontal eigenvalues, the vertical pairing is nonzero,
so the orthogonality relation yields a corresponding horizontal
identity. Taking finite linear combinations in these identities
allows us to use the spectral projections described above.
For any prescribed index shift $\ell\in\mathbb{Z}^2$, the
projections of plane waves with indices $n+\ell$ and $n$ recover
the Fourier mode $e^{i\ell\cdot\tilde{x}}$ as $|n|\to\infty$.
Every Fourier coefficient of $q_1-q_2$ therefore vanishes,
establishing uniqueness. The result concerns indices depending on
both horizontal variables at a single fixed wavenumber and a
fixed quasiperiodicity parameter; the high-energy limit involves
only auxiliary horizontal spectral parameters and does not
require additional measurement frequencies. 

The paper is organized as follows. Section 2 introduces the problem formulation and the geometric model. In Section 3 we derive a fundamental orthogonality relation for the difference $q_{1}-q_{2}$ from the near-field data.  Section 4 introduces the
associated quasiperiodic Schr\"odinger operators and establishes
the plane-wave approximation by finite spectral projections. 
In Section 5 we combine the orthogonality relation  with the
spectral approximation  obtained in Section 4 and conclude that it forces $q_{1}=q_{2}$ almost everywhere in the periodic cell, thus establishing our main uniqueness theorem for the inverse problem.

\section{Problem Formulation}
	
	We consider time-harmonic acoustic wave scattering by a bi-periodic layered medium occupying the upper half-space $\mathbb{R}_{+}^{3}:=\mathbb{R}^{2}\times(0,\infty)$. The lower boundary $\Gamma_{0} := \mathbb{R}^{2}\times\{ 0\}$ is assumed to be a sound-soft plane. We suppose that the medium is homogeneous and isotropic in the region above the plane $\{x_{3}=b\}$. The incident wave is supposed to be a time-harmonic wave of the form $u^{in}(x)e^{-i\omega t}$ with angular frequency $\omega>0$. The incident field $u^{in}$ impinges on $\Gamma_{b}$ from above and takes the plane wave form:
	\[
	u^{in}(x)=e^{ikx\cdot \hat{\theta}},\quad i=\sqrt{-1},
	\]
	with the incident direction $\hat{\theta}=(\sin\theta_{1}\cos\theta_{2},\sin\theta_{1}\sin\theta_{2},-\cos\theta_{1})$, where the incidence angles $\theta_{1}\in(-\frac{\pi}{2},\frac{\pi}{2}),\theta_{2}\in(0,2\pi]$. Here, $k:=\frac{\omega}{c_{0}}$ denotes the background wavenumber with $c_{0}>0$ the speed of sound in the homogeneous background medium. We can rewrite the incident field in the form $u^{in}(x)=e^{i(\alpha\cdot\tilde{x}-\beta x_{3})}$ with $\alpha=k\sin\theta_{1}(\cos\theta_{2},\sin\theta_{2}),\beta=k\cos\theta_{1}$. The inhomogeneous bi-periodic medium can be characterized by the  refractive index $q\in L^{\infty}(\mathbb{R}^{3})$, which is assumed to be $2\pi$-periodic in $\tilde{x}=(x_{1},x_{2})$ and to satisfy $q\equiv1$ in $x_{3}>b$.
	
	The scattering problem of the total field $u$ can be modeled by
	\begin{align}
		& \Delta u + k^2q u = 0, && \text{in } \mathbb{R}^3_+, \label{1.1} \\
		& u = 0, && \text{on } \Gamma_0, \label{1.3} \\
		& u = u^{in} + u^{sc}, && \text{in } \mathbb{R}^3_+. \label{1.4}
	\end{align}
	
	A function $v$ is called $\alpha$-quasiperiodic in $\tilde{x}$ if $e^{-i\alpha\tilde{x}}v(x)$ is $2\pi$-periodic in $\tilde{x}$. Equivalently,
	\begin{equation*}
		v(x_{1}+2\pi n_{1},x_{2}+2\pi n_{2} ,x_{3})=e^{i2\pi\alpha\cdot n}v(x_{1},x_{2},x_{3}),\quad n=(n_{1},n_{2})\in\mathbb{Z}^{2}.
		\label{ALPHA QUASIPERIODIC}
	\end{equation*}
	Since the structure is periodic in the $\tilde{x}$-direction and the incident wave $u^{in}$ is $\alpha$-quasiperiodic, we require the same for $u^{sc}$. Then, the function $u^{sc}e^{-i\alpha\cdot\tilde{x}}$ is $2\pi$-periodic with respect to $\tilde{x}$, and hence can be expanded into the Fourier series:
	\[
	u^{sc}(x)\,e^{-i\alpha\cdot\tilde{x}}=\sum\limits_{n\in\mathbb{Z}^{2}}u_{n}(x_{3})e^{in\cdot\tilde{x}},\quad x_{3}>0,
	\]
	which gives
	\begin{equation}
		u^{sc}(x)=\sum\limits_{n\in\mathbb{Z}^{2}}u_{n}(x_{3})e^{i\alpha_{n}\cdot\tilde{x}},\quad\alpha_{n}=n+\alpha\in\mathbb{R}^{2}.
		\label{scattered filed}
	\end{equation}
	In the homogeneous region, the scattered field $u^{sc}$ satisfies the Helmholtz equation
	\begin{equation}
		\Delta u^{sc}+k^{2}u^{sc}=0 \quad\text{in}\quad x_{3}>b.
	\end{equation}
	Substituting the Fourier expansion (\ref{scattered filed}) into this equation, each Fourier coefficient $u_{n}(x_{3})$ must satisfy the one-dimensional ordinary differential equation:
	\begin{equation}
		u_{n}^{''}(x_{3})+(k^{2}-|\alpha_{n}|^{2})u_{n}(x_{3})=0,\quad x_{3}>b.
		\label{ODE}
	\end{equation}
	Define the propagation constants:
	\begin{equation*}
		\beta_{n}:=\left\{\begin{array}{ll}
			\sqrt{k^{2}-|\alpha_{n}|^{2}} &\text{if}\quad \left|\alpha_{n}\right|<k ,\\
			i \sqrt{|\alpha_{n}|^{2}-k^{2}} &\text{if}\quad  \left|\alpha_{n}\right| \geq k.
		\end{array}\right.
	\end{equation*}
	Throughout we assume that Wood anomalies are excluded, i.e. $\beta_{n}\ne 0$ for all $n\in\mathbb{Z}^{2}$. The general solution of (\ref{ODE}) is as follows:
	\begin{equation}
		u_{n}(x_{3})=A_{n}e^{i\beta_{n}x_{3}}+B_{n}e^{-i\beta_{n}x_{3}},\quad A_{n},B_{n}\in\mathbb{C}.   \label{general solution}
	\end{equation}
	Combining with (\ref{scattered filed}) and (\ref{general solution}), we get the expression of the scattered field as a superposition of plane waves:
	\[
	u^{sc}(x)=\sum\limits_{n\in\mathbb{Z}^{2}}A_{n}e^{i\alpha_{n}\tilde{x}+i\beta_{n}x_{3}}+B_{n}e^{i\alpha_{n}\tilde{x}-i\beta_{n}x_{3}}.
	\]
	From a physical point of view, the scattered field $u^{sc}$ must remain bounded as $x_3\to\infty$. Thus, $u^{sc}$ consists only of outgoing or exponentially decaying waves, which leads to the known Rayleigh expansion condition,
	\begin{equation}
		u^{sc}(x)=\sum_{n\in \mathbb{Z}^{2}} u_{n}\, e^{i \alpha_{n} \cdot\tilde{x}+i \beta_{n} (x_{3}-b)}, \quad x_{3}>b,
		\label{rayleigh expansion1}
	\end{equation}
	where $u_{n}\in\mathbb{C}$ are the Rayleigh coefficients. The previous expansion \eqref{rayleigh expansion1} will serve as the upward radiation condition of the scattered field in this paper. For notational convenience, we write $\alpha=\alpha_{0}\in\mathbb{R}^{2},\beta=\beta_{0}\in\mathbb{C}$. The total field can then be expressed as
	\begin{equation*}
		u(x)=e^{i\alpha\cdot\tilde{x}-i\beta x_{3}}+\sum_{n\in \mathbb{Z}^{2}} u_{n} e^{i \alpha_{n} \cdot\tilde{x}+i \beta_{n} (x_{3}-b)},\quad x_{3}>b.
	\end{equation*}
	
		The direct scattering problem is to determine the solution of problem (\ref{1.1})-(\ref{1.4}) together with the radiation condition (\ref{rayleigh expansion1}). This can be achieved through a variational formulation with the Dirichlet-to-Neumann (DtN) opertator in appropriate Sobolev spaces. To this end, we introduce
		the relevant quasi-periodic and periodic Sobolev spaces and the DtN operator.  Setting
	
	\[
	\Omega_{b}:=(0,2\pi)^{2}\times(0,b),\quad\Omega_{b}^{c}=(0,2\pi)^{2}\times(b,\infty),\quad\Gamma_{b}=(0,2\pi)^{2}\times\{b\},
	\]
	we define
	\begin{equation*}
		\begin{aligned}
			H_{\alpha}^{1}(\Omega_{b})&:=\{ u(x)\in H^{1}(\Omega_{b}):e^{-i\alpha\cdot\tilde{x}}u(x) \text{ is }2\pi\text{-periodic w.r.t. }\tilde{x}\},\\
			H_{\alpha,0}^{1}(\Omega_{b})&=\{u(x)\in H_{\alpha}^{1}(\Omega_{b}),u=0 \text{ on }  \Gamma_{0}\},\\
			H_{\alpha,0,loc}^{1}(\mathbb{R}_{+}^{3})&=\{u\in H_{loc}^{1}(\mathbb{R}_{+}^{3}): u \text{ is $\alpha$-quasiperiodic w.r.t $\tilde{x}$}, u=0 \text{ on }\Gamma_{0}  \},\\
			H_{p}^{1}(\Omega_{b})&:=\{ v(x)\in H^{1}(\Omega_{b}):v(x) \text{ is }2\pi\text{-periodic w.r.t. }\tilde{x}\},\\
			H_{p,0}^{1}(\Omega_{b})&=\{v(x)\in H_{p}^{1}(\Omega_{b}),v=0 \text{ on }  \Gamma_{0}\}.\\
		\end{aligned}
	\end{equation*}
	For $s\in\mathbb{R}$, define the quasi-periodic Sobolev space on $\Gamma_{b}$ by
	\[
	H_{\alpha}^{s}(\Gamma_{b}):=\{f(\tilde{x})=\sum\limits_{n\in\mathbb{Z}^{2}}f_{n}e^{i\alpha_{n}\cdot\tilde{x}}:\sum\limits_{n\in\mathbb{Z}^{2}}|f_{n}|^{2}(1+|\alpha_{n}|^{2})^{s}<\infty\}
	\]
	equipped with the norm
	\[
	||f||_{H_{\alpha}^{s}(\Gamma_{b})}^{2}:=4\pi^{2}\sum\limits_{n\in\mathbb{Z}^{2}}|f_{n}|^{2}(1+|\alpha_{n}|^{2})^{s},
	\]
	where $f_{n}$ are the Fourier coefficients of the function $f$.
	In particular, the trace space of
	$H^1_\alpha(\Omega_b)$ on $\Gamma_b$ is given by $H^{1/2}_\alpha(\Gamma_b)$.

	For $\phi\in H_{\alpha}^{\frac{1}{2}}(\Gamma_{b})$, define the exterior DtN operator $\mathbb{T}:H_{\alpha}^{\frac{1}{2}}(\Gamma_{b})\to H_{\alpha}^{-\frac{1}{2}}(\Gamma_{b})$ by
	\[
	(\mathbb{T}\phi)(\tilde{x}):=\sum\limits_{n\in\mathbb{Z}^{2}}i\beta_{n}\phi_{n}e^{i\alpha_{n}\cdot\tilde{x}},\quad \phi(\tilde{x})=\sum\limits_{n\in\mathbb{Z}^{2}}\phi_{n}e^{i\alpha_{n}\cdot\tilde{x}},
	\]
	where $\phi_{n}=\frac{1}{4\pi^{2}}\int_{0}^{2\pi}\int_{0}^{2\pi}\phi(\tilde{x})e^{-i\alpha_{n}\cdot\tilde{x}}dx_{1}dx_{2}$. The expression of $\mathbb{T}$ stems from the Neumann trace of $\frac{\partial w}{\partial x_{3}}$ on $\Gamma_{b}$, where $w\in H_{\alpha,loc}^{1}(\Omega_{b}^{c})$ is the unique solution of the following exterior problem:
	\begin{equation*}
		\begin{cases}
			\Delta w+k^{2} w=0 & \quad\text{in }\Omega_{b}^{c}, \\
			w=\phi & \quad\text{on }\Gamma_{b},\\
			w=\sum\limits_{n\in\mathbb{Z}^{2}}w_{n}\,e^{i\alpha_{n}\cdot\tilde{x}+i\beta_{n}(x_{3}-b)}& \quad\text{in }x_{3}>b.
		\end{cases}
	\end{equation*}
	
	From the asymptotics $\beta_{n}=i|\alpha_{n}|+O(\frac{1}{|n|}), |n|\to\infty$, the operator $\mathbb{T}:H_{\alpha}^{\frac{1}{2}}(\Gamma_{b})\to H_{\alpha}^{-\frac{1}{2}}(\Gamma_{b})$ is well defined and bounded. Since the scattered field $u^{sc}=u-u^{in}$ satisfies the Rayleigh expansion (\ref{rayleigh expansion1}), we have
	\[
	\mathbb{T}\left( \left. u^{sc}(x) \right|_{\Gamma_b} \right)
	=\left. \frac{\partial u^{sc}}{\partial x_3} \right|_{\Gamma_b}.
	\]
	We further get the transparent boundary condition for the total field $u$ on $\Gamma_{b}$:
	\begin{equation*}
		\begin{aligned}
			\left. \frac{\partial u(x)}{\partial x_{3}}\right|_{\Gamma_b}&=\left. \frac{\partial u^{sc}(x)+\partial u^{in}(x)}{\partial x_{3}}\right|_{\Gamma_b}=\mathbb{T}\left( \left. (u(x) - u^{in}(x)) \right|_{\Gamma_b} \right)+\left. \frac{\partial u^{in}(x)}{\partial x_{3}}\right|_{\Gamma_b}\\
			&=\mathbb{T}\left( \left. u \right|_{\Gamma_b} \right)-2i\beta e ^{i\alpha\cdot\tilde{x}}e^{-i\beta b}.
		\end{aligned}
	\end{equation*}
	The direct scattering problem can be equivalently formulated as the following boundary value problem over the truncated domain $\Omega_{b}$:
	\begin{equation}
		\begin{cases}
			\Delta u+k^{2} q(x) u=0 & \quad\text{in }\Omega_{b} \\
			u=0 & \quad\text{on }\Gamma_{0}, \\
			\frac{\partial u}{\partial x_{3}}=\mathbb{T}u+f & \quad\text{on }\Gamma_{b},
		\end{cases}
		\label{direct problem.}
	\end{equation}
	where $f=-2i\beta e ^{i\alpha\cdot\tilde{x}}e^{-i\beta b}$.  Multiplying the equation (\ref{direct problem.}) by the complex conjugate of a test function $v\in H_{\alpha,0}^{1}(\Omega_{b})$, integrating over $\Omega_{b}$ and using Green's identity, we obtain the following variational problem: find $u\in H_{\alpha,0}^{1}(\Omega_{b})$ such that
	\begin{equation}
		a(u,\varphi)=(f,\varphi), \quad \text{for all } \varphi\in H_{\alpha,0}^{1}(\Omega_{b}),
		\label{variation problrm}
	\end{equation}
	with the sesquilinear form
	\begin{equation}
		a(u, \varphi):=\int_{\Omega_{b}}\nabla u\cdot\nabla\overline{\varphi}-k^{2}qu\overline{\varphi}dx-\int_{\Gamma_{b}}\mathbb{T}u\cdot\overline{\varphi}ds,
		\label{SES FORM}
	\end{equation}
	and the bounded linear functional
	\[
(f,\varphi)=-2i\beta e^{-i\beta b}\int_{\Gamma_{b}}e^{i\alpha\cdot\tilde{x}}\overline{\varphi}ds.
	\]
	Based on the above variational formulation, we collect the existence and uniqueness results to the direct scattering problem as follows.
	\begin{proposition}
		Let $\alpha=\alpha(k)=k\sin\theta_{1}(\cos\theta_{2},\sin\theta_{2})$, where $\theta_{1}\in(-\frac{\pi}{2},\frac{\pi}{2}),\theta_{2}\in(0,2\pi]$.\\
		(i) For any $k>0$, the direct scattering problem (\ref{1.1})-(\ref{1.4})with the incident plane wave $u^{\mathrm{in}}(x)=e^{i\alpha\cdot\tilde{x}-i\beta x_3}$ has at least one solution $u\in H_{\alpha,0,loc}^{1}(\mathbb{R}_{+}^{3})$ whose scattered field  satisfies the Rayleigh radiation condition (\ref{rayleigh expansion1}).\\
		(ii) There exists a unique solution if one of the following conditions holds true:
		\begin{enumerate}
			\item $k\in\mathbb{R}_{+}\backslash\mathcal{D}$, where $\mathcal{D}\subset\mathbb{R}$ is a discrete set whose only accumulation point is at infinity;
			\label{con1}
			\item $q\in W^{1,\infty}(\Omega_b;\mathbb R)$ is real valued and there exists $c\ge0$ such that $(x_{3}+c)(1-q(x))\ge0$ on $\Gamma_{b}$ and $(x_{3}+c)\frac{\partial q}{\partial x_{3}}\ge0$ a.e. in $\Omega_{b}$.
			\label{con2}
			\item ${\rm Re }\, q\geq 0$  in $\mathbb{R}^3_+$ and ${\rm Im }\, q\ge0$ in $\Omega_{b}$. There exists $q_{0}$ such that  ${\rm Im}\, q(x)\geq q_0>0$ in $B_\epsilon(x_0)$ for some $x_0\in \Omega_b$ and some $\epsilon>0$.
			\item  \(k_\varepsilon=k+i\varepsilon\) with \(k>0\) and \(\varepsilon>0\) and   $q(x)\geq \sin^2\theta_1$ in $\mathbb{R}^3_+$. That is, the  complex-frequency scattering problem is uniquely solvable.
			\item $q(x)\ge\sin^{2}\theta_{1}$ and
			the total field satisfies the additional orthogonal relation:
			\begin{equation*}
			\int_{Q_{\infty}^{+}}\left [ \alpha\cdot\nabla_{\tilde{x}}u-ik^{2}q(x)u \right ] \overline{\phi}dx=0,\quad Q_{\infty}^{+}:=(0,2\pi)^2\times(0,\infty),
			\end{equation*}\label{con3}
			for all $\phi\in H_{\alpha,0,loc}^{1}(\mathbb{R}_{+}^{3})$ solving the homogeneous scattering problem with $u^{in}=0$. This solution is the limiting-absorption solution associated with assertion 4.
		
		Under any of conditions1--5, problem (\ref{direct problem.}) is uniquely solvable for every \(f\in H^{-1/2}_\alpha(\Gamma_b)\), with
		\[
		\|u\|_{H_{\alpha}^1(\Omega_b)}
		\le C\|f\|_{H^{-1/2}_\alpha(\Gamma_b)},
		\]
		where \(C\)  is a constant independent of \(f\). Assertion 5 gives uniqueness within the class satisfying the additional constraint and does not exclude nontrivial homogeneous outgoing solutions.
			
		\end{enumerate}
		\label{prop2.1}
	\end{proposition}
The first assertion of Proposition \ref{prop2.1}(i) shows existence of solutions with plane-wave incidence, while each condition of the second assertion  (ii) gives uniqueness and existence, that is, the direct problem (\ref{1.1})-(\ref{1.4}) admits a unique solution, which can be obtained by applying analytic Fredholm theory. More precisely, for fixed $k>0$ and $\alpha\in\mathbb{R}^{2}$, the sesquilinear form (\ref{SES FORM}) defines a Fredholm operator of index zero. By the Fredholm alternative, if the homogeneous variational problem $a(u,\varphi)=0$ has only a trivial solution $u=0$ for all $\varphi\in H_{\alpha,0}^{1}(\Omega_{b})$, then the variational problem (\ref{variation problrm}) admits a unique solution $u\in H_{\alpha,0}^{1}(\Omega_{b})$. Moreover, the solution continuously depends on the incident field. For fixed incident angles, if $\alpha=\alpha(k)$ varies with wavenumbers, then the function space varies with $k$. We need to transform the quasi-periodic space into a periodic space to ensure that the function space is independent of the wave number $k$. Let $v(x)=e^{-i\alpha\cdot\tilde{x}}u(x)$, then $v(x)\in H_{p,0}^{1}(\Omega_{b})$. By the Riesz representation theorem, the variational form (\ref{variation problrm}) can be rewritten as the operator equation
	\[
	(\mathcal{T}-\mathcal{K}_{k})v=F_{k},
	\]
	where the operator $\mathcal{T}$ is coercive. The compact operator $\mathcal{K}_{k}:H_{p,0}^{1}(\Omega_{b})\to H_{p,0}^{1}(\Omega_{b})^{\prime}$ depends analytically on $k\in\mathbb{C}$ such that $|\alpha_{n}|\ne k$ for all $n\in\mathbb{Z}^{2}$. Here $H_{p,0}^{1}(\Omega_{b})^{\prime}$ denotes the dual space of $H_{p,0}^{1}(\Omega_{b})$. Together with the analytic Fredholm theory, we can prove that the operator $\mathcal{T}-\mathcal{K}_{k}$ is invertible for all $k\in\mathbb{C}\backslash \mathcal{D}$, where $\mathcal{D}$ is a discrete set. We therefore obtain that the scattering problem (\ref{1.1})-(\ref{1.4}) admits a unique solution for all but a discrete set of wave numbers $k$. Equivalently, the variational problem (\ref{variation problrm}) is uniquely solvable in $H_{\alpha,0}^{1}(\Omega_{b})$. For more details we refer to \cite{BL2025,BBS1994,ElschnerSchmidt1998,Kirsch A1995,Strycharz B 1998}.
	
	 For a real-valued refractive index, taking the imaginary part of the variational problem (\ref{variation problrm}) with $u^{in}=0$ yields $u_n=0$ for all $|\alpha_n|<k$, that is, all propagating Rayleigh amplitudes $u_{n}$ must vanish. Hence, any non-trivial homogeneous outgoing solution, if it exists, can contain only evanescent Rayleigh modes. This process together with the unique continuation argument also proves uniqueness under the condition (ii) (3).
	The condition in (ii) (2) is another sufficient condition on the refractive index for uniqueness.
	 The Rellich identity associated with condition (ii) (2) then forces such a solution to vanish identically. Equivalently, the corresponding homogeneous scattering operator has a trivial kernel, and hence yields uniqueness of the direct scattering problem \cite{BBS1994,CMT, ElschnerSchmidt1998, HuLiuQuZhang2015,Kirsch A1995,Strycharz B.,Zhang B2018}. By contrast, the orthogonality condition in (ii) (5) arises from the limiting absorption principle. More precisely, the wavenumber $k$ is replaced by $k_\varepsilon=k+i\varepsilon$ with $\varepsilon>0$, for which the corresponding dissipative scattering problem is uniquely solvable \cite{HuKi26,HuKiZh,Kirsch A2018}, which proves the sufficient condition (ii) (4). At an exceptional real wave number, the homogeneous limiting operator has a non-trivial kernel. Passing to the limit $\varepsilon\to0^+$, the solvability of the dissipative scattering problem requires the first-order perturbation term to be orthogonal to this kernel, which yields condition (\ref{con3}). Hence, (\ref{con3}) serves as the additional compatibility condition that characterizes the limiting-absorption solution at frequencies where the Rayleigh radiation condition alone does not ensure uniqueness. Thus, each condition in the second asserion (ii) excludes guided modes, whereas condition (\ref{con3}) ensures uniqueness when non-trivial guided modes may be present.
	
	We choose $k\notin \mathcal{D}$ or suppose that $q$ or $u$ satisfies one of the conditions in Proposition \ref{prop2.1} to ensure that the direct scattering problem is uniquely solvable. In our inverse problem, we additionally assume that the refractive index  $q$ is real-valued and  depends only on $\tilde{x}$.
	
%

	Below we consider an $\alpha$-quasiperiodic point source wave $G(x,y)$ for the inverse problem. Let $y\in\Gamma_{a}:=\{(\tilde{x},a): 0<x_{1},x_{2}<2\pi\}$ with some fixed $a>b$. 
	 The $\alpha$-quasiperiodic Green function $G(x,y)$ in the homogeneous medium is given by
	\begin{equation}
		G(x,y)=\frac{i}{8\pi^{2}}\sum_{n\in\mathbb{Z}^{2}}\frac{1}{\beta_{n}}\exp(i\alpha_{n}\cdot(\tilde{x}-\tilde{y})+i\beta_{n}|x_{3}-y_{3}|),\quad x-y\notin2\pi\mathbb{Z}^{2}\times\left\{0\right\}.
		\label{G(x,y)}
	\end{equation}
	
	The well-posedness result for general boundary data also applies
	to this incident field. Indeed, since $y_3=a>b$, its boundary
	datum is
	\[
	f_y
	:=\left.(\partial_{x_3}-\mathbb{T})G(\cdot,y)\right|_{\Gamma_b}
	=\frac{1}{4\pi^2}\sum_{n\in\mathbb Z^2}
	\exp\!\left(
	i\alpha_n\cdot(\tilde{x}-\tilde{y})
	+i\beta_n(a-b)
	\right).
	\]
	There are only finitely many propagating modes, while the
	coefficients of the evanescent modes decay exponentially
	as $|n|\to\infty$ since $a>b$. Hence
	$f_y\in H^{-1/2}_{\alpha}(\Gamma_b)$. Applying the general
	boundary-data result in Proposition 2.1 with $f=f_y$ gives
	a unique solution $u(\cdot,y)\in H^1_{\alpha,0}(\Omega_b)$,
	satisfying
	\[
	\|u(\cdot,y)\|_{H^1_{\alpha}(\Omega_b)}
	\le C\|f\textit{}_y\|_{H^{-1/2}_{\alpha}(\Gamma_b)},
	\]
	where $C$ is independent of $y$. The outgoing extension of
	$u(\cdot,y)|_{\Gamma_b}-G(\cdot,y)|_{\Gamma_b}$ defines the
	scattered field above $\Gamma_b$, and hence yields the
	 scattering solution corresponding to a  point-source incidence.
	
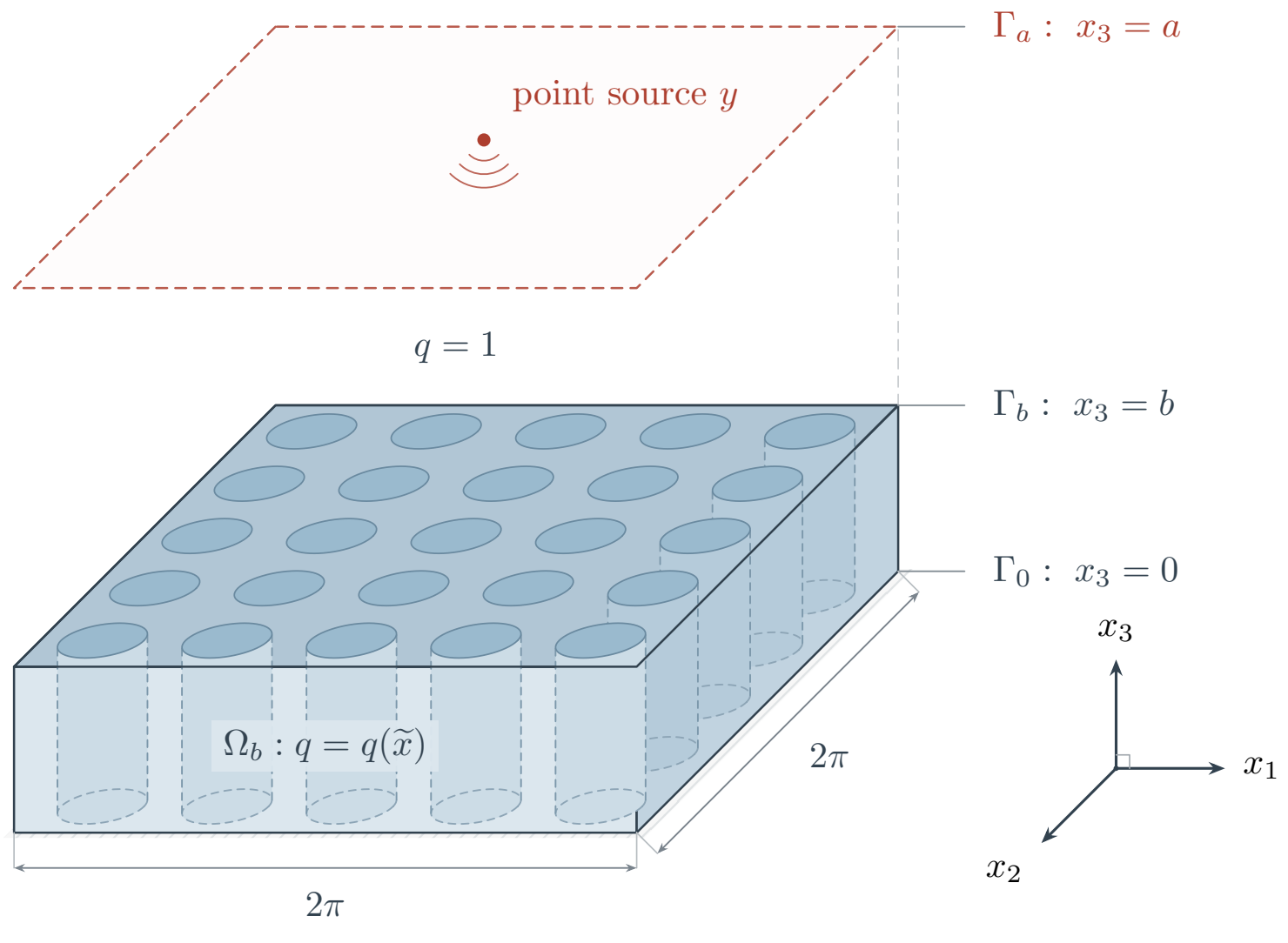
\begin{figure}[htbp]
		\centering
		\newcommand{\MaxwellPlotPlanarFinal}{%
\begin{tikzpicture}[
	x={(1cm,0cm)}, y={(-0.42cm,-0.42cm)}, z={(0cm,1cm)},
	font=\small,line cap=round,line join=round,
	outline/.style={draw=ink,line width=0.60pt},
	guide/.style={draw=ink!28,line width=0.4pt,densely dashed},
	dimension/.style={draw=ink!65,line width=0.45pt,
		{Stealth[length=3pt]}-{Stealth[length=3pt]}},
	axis/.style={draw=ink,line width=0.75pt,
		-{Stealth[length=4.5pt,width=3.5pt]}},
	leader/.style={draw=ink!60,line width=0.45pt}
	]
	\definecolor{ink}{RGB}{50,66,79}
	\definecolor{topblue}{RGB}{177,198,213}
	\definecolor{leftblue}{RGB}{220,231,238}
	\definecolor{rightblue}{RGB}{193,211,224}
	\definecolor{profileblue}{RGB}{107,153,183}
	\definecolor{profileedge}{RGB}{80,117,143}
	\definecolor{accent}{RGB}{173,62,46}
	\def\L{6}
	\def\B{1.60}
	\def\A{5.25}
	\def\radius{0.40}
	\coordinate (O) at (0,0,0);
	\coordinate (F) at (0,\L,0);
	\coordinate (C) at (\L,\L,0);
	\coordinate (R) at (\L,0,0);
	\coordinate (Ob) at (0,0,\B);
	\coordinate (Fb) at (0,\L,\B);
	\coordinate (Cb) at (\L,\L,\B);
	\coordinate (Rb) at (\L,0,\B);
	\coordinate (Oa) at (0,0,\A);
	\coordinate (Fa) at (0,\L,\A);
	\coordinate (Ca) at (\L,\L,\A);
	\coordinate (Ra) at (\L,0,\A);
	
	\path[fill=black!3]
	(-0.05,-0.05,0)--(6.12,-0.05,0)--(6.12,6.12,0)--(-0.05,6.12,0)--cycle;
	\path[pattern=north east lines,pattern color=ink!15]
	(-0.05,-0.05,0)--(6.12,-0.05,0)--(6.12,6.12,0)--(-0.05,6.12,0)--cycle;
	
	\path[fill=leftblue,outline] (F)--(C)--(Cb)--(Fb)--cycle;
	\path[fill=rightblue,outline] (C)--(R)--(Rb)--(Cb)--cycle;
	\path[fill=topblue,outline] (Ob)--(Fb)--(Cb)--(Rb)--cycle;
	
	\pgfmathsetmacro{\thetaR}{atan(-0.42)}
	\pgfmathsetmacro{\thetaL}{\thetaR+180}
	\foreach \cx/\cy in {5.40/0.60,5.40/1.80,5.40/3.00,5.40/4.20,
		0.60/5.40,1.80/5.40,3.00/5.40,4.20/5.40,5.40/5.40}{
		\pgfmathtruncatemacro{\frontrow}{abs(\cy-5.4)<0.001}
		\ifnum\frontrow=1
		\colorlet{profilebody}{profileblue!14!leftblue}
		\else
		\colorlet{profilebody}{profileblue!14!rightblue}
		\fi
		\path[fill=profilebody]
		plot[domain=\thetaR:\thetaL,samples=31,variable=\ang]
		({\cx+\radius*cos(\ang)},{\cy+\radius*sin(\ang)},\B)
		--({\cx+\radius*cos(\thetaL)},{\cy+\radius*sin(\thetaL)},0)
		--plot[domain=\thetaL:\thetaR,samples=31,variable=\ang]
		({\cx+\radius*cos(\ang)},{\cy+\radius*sin(\ang)},0)
		--cycle;
		\path[fill=profilebody]
		plot[domain=0:360,samples=57,variable=\ang]
		({\cx+\radius*cos(\ang)},{\cy+\radius*sin(\ang)},0)--cycle;
		\draw[draw=profileedge!65,line width=0.45pt,
		dash pattern=on 2pt off 1.6pt]
		plot[domain=0:360,samples=57,variable=\ang]
		({\cx+\radius*cos(\ang)},{\cy+\radius*sin(\ang)},0);
		\foreach \th in {\thetaR,\thetaL}{
			\draw[draw=profileedge!75,line width=0.45pt,
			dash pattern=on 2pt off 1.6pt]
			({\cx+\radius*cos(\th)},{\cy+\radius*sin(\th)},0)
			--({\cx+\radius*cos(\th)},{\cy+\radius*sin(\th)},\B);
		}
	}
	
	\foreach \uu in {0.60,1.80,3.00,4.20,5.40}{
		\foreach \vv in {0.60,1.80,3.00,4.20,5.40}{
			\path[fill=profileblue!68!white,draw=profileedge!85,line width=0.42pt]
			plot[domain=0:360,samples=57,variable=\ang]
			({\uu+\radius*cos(\ang)},{\vv+\radius*sin(\ang)},\B)
			--cycle;
		}
	}
	\draw[outline] (Ob)--(Fb)--(Cb)--(Rb)--cycle;
	
	\path[fill=accent,fill opacity=0.015] (Oa)--(Fa)--(Ca)--(Ra)--cycle;
	\draw[draw=accent!85,line width=0.60pt,dash pattern=on 3pt off 2.2pt]
	(Oa)--(Fa)--(Ca)--(Ra)--cycle;
	\draw[guide] (Rb)--(Ra);
	\node[text=ink] at (3,3,3.42) {$q=1$};
	
	\coordinate (S) at (3.1,2.6,\A);
	\fill[accent] (S) circle[radius=1.8pt];
	\node[above right=4pt,text=accent] at (S) {point source $y$};
	\begin{scope}[shift={(S)},x={(1cm,0cm)},y={(0cm,1cm)}]
		\foreach \rr in {0.20,0.33,0.46}{
			\draw[draw=accent!80,line width=0.45pt]
			(225:\rr) arc[start angle=225,end angle=315,radius=\rr cm];
		}
	\end{scope}
	
	\foreach \pt in {R,Rb,Ra}{
		\draw[leader] (\pt)--([xshift=0.64cm]\pt);
	}
	\node[anchor=west,xshift=0.78cm,text=ink] at (R)
	{$\Gamma_0:\ x_3=0$};
	\node[anchor=west,xshift=0.78cm,text=ink] at (Rb)
	{$\Gamma_b:\ x_3=b$};
	\node[anchor=west,xshift=0.78cm,text=accent] at (Ra)
	{$\Gamma_a:\ x_3=a$};
	\node[text=ink,fill=leftblue,fill opacity=0.96,text opacity=1,
	inner xsep=3pt,inner ysep=2pt] at (3,6,0.84)
	{$\Omega_b:q=q(\widetilde x)$};
	
	\coordinate (D1) at ([yshift=-0.34cm]F);
	\coordinate (D2) at ([yshift=-0.34cm]C);
	\coordinate (D3) at ([xshift=0.20cm,yshift=-0.20cm]C);
	\coordinate (D4) at ([xshift=0.20cm,yshift=-0.20cm]R);
	\draw[draw=ink!35,line width=0.4pt] (F)--(D1) (C)--(D2) (C)--(D3) (R)--(D4);
	\draw[dimension] (D1)--(D2)
	node[midway,below=3pt,text=ink] {$2\pi$};
	\draw[dimension] (D3)--(D4)
	node[midway,below right=2pt,text=ink] {$2\pi$};
	
	\coordinate (Axes) at ([xshift=2.1cm,yshift=-1.9cm]R);
	\begin{scope}[shift={(Axes)},x={(1cm,0cm)},y={(0cm,1cm)}]
		\draw[axis] (0,0)--(-0.72,-0.72) node[below left=1pt] {$x_2$};
		\draw[axis] (0,0)--(1.05,0) node[right=1pt] {$x_1$};
		\draw[axis] (0,0)--(0,1.05) node[above=1pt] {$x_3$};
		\draw[draw=ink!45,line width=0.4pt] (0,0.13)--(0.13,0.13)--(0.13,0);
		\fill[ink] (0,0) circle[radius=0.7pt];
	\end{scope}
\end{tikzpicture}
}
\resizebox{0.6\linewidth}{!}{\MaxwellPlotPlanarFinal}
\caption{Geometric model of the grating diffraction problem (one periodic cell).}
\label{fig:refined_grating_2d_model}
\end{figure}

	The inverse problem we consider is to determine the bi-periodic refractive index function $q$ from the near-field data $u(x,y)$ for all $x, y\in\Gamma_{a}$ with $x\neq y$, where $u(\cdot,y)$ denotes the total field of the problem (\ref{1.1})-(\ref{1.4}) subject to the Rayleigh expansion radiation condition (\ref{rayleigh expansion1}) and with the given incident wave $u^{in}(x,y)=G(x,y)$. Throughout this paper we focus on the unique identification problem of this inverse problem; see Figure~\ref{fig:refined_grating_2d_model}.

	\section{The orthogonality relation}
	
	In this section, we derive an orthogonality identity for the potential function $q$. First, we introduce an auxiliary Robin boundary value problem. Then, we establish a density result regarding the boundary traces of the total fields generated by point sources on the measurement plane $\Gamma_{a}$. This density result is then used to show that the coincidence of the near-field generated by two refractive indexes $q_{j} (j=1,2)$ implies that the interior Robin to Dirichlet(RtD) maps  coincide on the artificial boundary $\Gamma_{b}$. Finally, applying Green's formula and the principle of unique continuation, we obtain the desired orthogonality relation.
	
	Fix $\eta>0$ and set
	$
	R_\eta w:=\partial_\nu w-i\eta w\in H_{\sigma,0}^{-\frac{1}{2}}(\Gamma_{b})$ ,where $\sigma\in\left\{\alpha,-\alpha\right\}\subset\mathbb{R}^{2}$ and $\nu$ denotes the unit outward normal to $\Omega_{b}$.  For  fixed $\alpha$, we next prove the uniqueness result for the following Robin boundary value problem.
\begin{lemma}\label{uniquessrobin}
	Let $q$ be real-valued and bounded in $\Omega_b$. For every
	$f\in H^{-1/2}_\sigma(\Gamma_b)$, the problem
\begin{equation}\label{robinproblem}
\begin{cases}
	\Delta w+k^2qw=0&\text{in }\Omega_b,\\
	w=0&\text{on }\Gamma_0,\\
	R_\eta w=f&\text{on }\Gamma_b
\end{cases}
\end{equation}
has a unique solution $w\in H_{\sigma,0}^{1}(\Omega_{b})$, and
\[
\|w\|_{H^{1}(\Omega_{b})}\le C\|f\|_{H_{\sigma}^{-\frac{1}{2}}(\Gamma_{b})}.
\]
\end{lemma}
\begin{proof}
	Multiplying both side of (\ref{robinproblem}) by the test function $v\in H_{\sigma,0}^{1}(\Omega_{b})$, and applying  integration by parts, we obtain the variational formulation:
	find $w\in H_{\sigma,0}^{1}(\Omega_{b})$ such that
\begin{equation}\label{robinvariational}
	b(w,v)=(f,v)
\qquad\text{for all }v\in H_{\sigma,0}^{1}(\Omega_{b}),
\end{equation}
	where
	\[
	b(w,v)
	=\int_{\Omega_b}
	\bigl(\nabla w\cdot\nabla\overline v
	-k^2qw\overline v\bigr)\,dx
	-i\eta\int_{\Gamma_b}
	w\,\overline{ v}\,ds.
	\]
	The remaining boundary terms cancel due to  the quasiperiodicity and the condition $v=0$ on $\Gamma_{0}$.
The Cauchy--Schwarz inequality
	and the trace theorem imply that
	\[
	|b(w,v)|
	\le C\|w\|_{H^1(\Omega_b)}\|v\|_{H^1(\Omega_b)},
	\qquad w,v\in H_{\sigma,0}^{1}(\Omega_{b}).
	\]
	Set $\lambda=1+k^2\|q\|_{L^\infty(\Omega_b)}$ and define
	\[
	b_{1}(w,v)
	:=b(w,v)+\lambda\int_{\Omega_b}w\overline v\,dx.
	\]
	Since $q$ is real-valued, we have
	\[
	\begin{split}
		\operatorname{Re}b_{1}(w,w)
		&=\int_{\Omega_b}|\nabla w|^2\,dx
		+\int_{\Omega_b}(\lambda-k^2q)|w|^2\,dx\\
		&\ge \|\nabla w\|_{L^2(\Omega_b)}^2
		+\|w\|_{L^2(\Omega_b)}^2
		=\|w\|_{H^1(\Omega_b)}^2.
	\end{split}
	\]
	Thus $b_{1}$ is coercive on $H_{\sigma,0}^{1}(\Omega_{b})$. Let
	$A,B,M:H_{\sigma,0}^{1}(\Omega_{b})\to H_{\sigma,0}^{1}(\Omega_{b})'$ be bounded linear operators  defined by
	\[
	\begin{aligned}
		(Aw)(v)=b(w,v),\quad
		(Bw)(v)=b_{1}(w,v),\quad
		(Mw)(v)=\int_{\Omega_b}w\overline v\,dx.
	\end{aligned}
	\]
	Then the variational problem (\ref{robinvariational}) can be rewritten as follows:
	\[
	(Aw,v)=((B-\lambda M)w,v)=F_{f}(v),
	\]
	with the bounded linear functional
	\[
	F_{f}(v)=(f,v)_{\Gamma_{b}}\quad\text{for all }\quad v\in H_{\sigma,0}^{1}(\Omega_{b}).
	\]
	The Lax--Milgram theorem shows that $B$ is an isomorphism.
	Since $H_{\sigma,0}^{1}(\Omega_{b})$ is a closed subspace of $H^1(\Omega_b)$ and
	$\Omega_b$ is bounded, the embedding
	$H_{\sigma,0}^{1}(\Omega_{b})\hookrightarrow L^2(\Omega_b)$ is compact. Consequently,
	$M:H_{\sigma,0}^{1}(\Omega_{b})\to H_{\sigma,0}^{1}(\Omega_{b})'$ is compact. It follows from
$
	A=B\bigl(I-\lambda B^{-1}M\bigr)
$
	that $A$ is Fredholm of index zero.
	By Fredholm alternative theorem, the uniqueness implies the existence of a unique solution to (\ref{robinproblem}). The  uniqueness is equivalently to the fact that  $\text{Ker} A=0$. Taking $v=w$
	in the variational equation and then taking the imaginary part,
	we obtain
	\[
	-\eta\int_{\Gamma_b}|w|^2\,ds=0.
	\]
	Since $\eta>0$, we have $w=0$ on $\Gamma_{b}$, and the Robin condition also gives
	$\partial_\nu w=0$ there. Since both the Dirichlet and Neumann traces of $w$ vanish on $\Gamma_b$,
	extending $w$ locally by zero across $\Gamma_b$ and applying unique continuation for Schrödinger equations with bounded potentials yields $w=0$ in $\Omega_b$; see \cite{Jerison D1985}. 
	That is,  $A$  is  an
	isomorphism from $H_{\sigma,0}^{1}(\Omega_{b})$ onto $H_{\sigma,0}^{1}(\Omega_{b})'$. For each
	$f\in H^{-1/2}_\sigma(\Gamma_b)$, the unique variational
	solution is $w=A^{-1}F_f$ and 
	\[
		\|w\|_{H^1(\Omega_b)}
		\le \|A^{-1}\|_{\mathcal L(H_{\sigma,0}^{1}(\Omega_{b})',H_{\sigma,0}^{1}(\Omega_{b}))}
		\|F_f\|_{H_{\sigma,0}^{1}(\Omega_{b})'}\le C\|f\|_{H^{-1/2}_\sigma(\Gamma_b)},
	\]
	where $C$ is a constant  independent of $f$. This completes the proof.
%
\end{proof}
 
  Based on  Lemma \ref{uniquessrobin}, we therefore can establish the density result as follows.
	\begin{lemma}
		Let $u(x,y)$ be the solution of the scattering problem (\ref{1.1})-(\ref{1.4}) with the Rayleigh condition (\ref{rayleigh expansion1}) corresponding to the incident field $u^{in}(x,y)=G(x,y)$. The span $\{ R_{\eta}u(\cdot,y)|_{\Gamma_b} \colon y \in \Gamma_a \}$ is dense in $H_{\alpha}^{-\frac{1}{2}}(\Gamma_b)$.
		
		\begin{proof} We follow the approach for the two dimensional problem in 
			\cite{Kirsch1994}.
			Let $\varphi\in H_{\alpha}^{\frac{1}{2}}(\Gamma_{b})$ satisfy $\int_{\Gamma_{b}}R_{\eta}u(x,y)\overline{\varphi}dx=0$ for all $y\in\Gamma_{a}$. By duality pairing, it suffices to demonstrate that $\varphi \equiv 0$. Let $v(x)\in H_{-\alpha,loc}^{1}(\mathbb{R}_{+}^{3})$ be the unique solution to the following transmission problem:
			\begin{equation}
				\begin{cases}
					\Delta v + k^{2}v = 0 & \text{in } \Omega_{b}^{c}, \\
					\Delta v + k^{2}q(x)v = 0 & \text{in } \Omega_{b}, \\
					v_{+}-v_{-}=\overline{\varphi} & \text{on } \Gamma_{b},\\
					(\frac{\partial v_{+}}{\partial\nu}-i\eta v_{+})-(\frac{\partial v_{-}}{\partial\nu} -i\eta v_{-})=0& \text{on } \Gamma_{b},\\
					v = 0 & \text{on } \Gamma_{0},
				\end{cases}
				\label{transmission problem1}
			\end{equation}
			where $\nu$ denotes the unit outward normal on $\Gamma_{b}$. The subscripts $+$ and $-$ indicate the traces taken from $\Omega_{b}^{c}$ and $\Omega_{b}$, respectively. The solution $v(x)$ of (\ref{transmission problem1}) is assumed to be $-\alpha$-quasiperiodic with the Rayleigh expansion condition:
			\begin{equation}
				v(x)=\sum_{n\in\mathbb{Z}^{2}}v_{n}e^{i(\hat{\alpha}_{n}\cdot\tilde{x}+\hat{\beta}_{n}x_{3})},x_{3}\ge b,
				\label{rayleigh expansion2}
			\end{equation}
			where $\hat{\alpha}_{n}=-\alpha+n\in \mathbb{R}^{2}$,$v_{n}\in \mathbb{C}$ and
			\begin{equation*}
				\hat{\beta}_{n}=\left\{\begin{array}{ll}
					\sqrt{k^{2}-|\hat{\alpha}_{n}|^{2}} &\text{if}\quad \left|\hat{\alpha}_{n}\right|<k ,\\
					i \sqrt{|\hat{\alpha}_{n}|^{2}-k^{2}} &\text{if}\quad  \left|\hat{\alpha}_{n}\right| \geq k.
				\end{array}\right.
			\end{equation*}
			For any $\varphi\in H_{\alpha}^{\frac{1}{2}}(\Gamma_{b})$, we have
			\begin{equation*}
				\begin{aligned}
					&\int_{\Gamma_{b}}R_{\eta}u(x,y)\overline{\varphi}ds(x)\\
					&=\int_{\Gamma_{b}}(\partial_{\nu}u(x,y)-i\eta u(x,y))\overline{\varphi}ds(x)\\
					&=\int_{\Gamma_{b}}\partial_{\nu}u(x,y)(v_{+}-v_{-})-u(x,y)(\partial_{\nu}v_{+}-\partial_{\nu}v_{-})ds(x)\\
					&=\int_{\Gamma_{b}}\left(u(x,y)\frac{\partial v_{-}}{\partial\nu}-v_{-}\frac{\partial u(x,y)}{\partial\nu}\right)ds(x)+\int_{\Gamma_{b}}\left(v_{+}\frac{\partial u(x,y)}{\partial\nu}-u(x,y)\frac{\partial v_{+}}{\partial\nu}\right)ds(x).
				\end{aligned}
			\end{equation*}
		Since both $u(x,y)$ and $v(x)$ satisfy the same Helmholtz equation $\Delta w+k^{2}q(x)w=0$ in $\Omega_{b}$ and vanish on $\Gamma_{0}$, 
		  Green's identity applied in $\Omega_{b}$ yields
			\begin{equation}
				\int_{\Gamma_{b}}\left(\frac{\partial u(x,y)}{\partial\nu}v_{-}-\frac{\partial v_{-}}{\partial\nu}u(x,y)\right)ds(x)=0. 
				\label{3.3}
			\end{equation}
Furthermore, since $u^{s}$ and $v_{+}$ satisfy the outgoing Rayleigh expansions (\ref{rayleigh expansion1}) and (\ref{rayleigh expansion2}) with modes $\alpha_{n}$ and $\hat{\alpha}_{n}=-\alpha+n$, respectively, and since $\hat{\beta}_{-n}=\beta_{n}$, we get
			\begin{equation}
				\int_{\Gamma_{b}} \left( u^{s}(x,y)\frac{\partial v_{+}}{\partial\nu} - v_{+}\frac{\partial u^{s}(x,y)}{\partial\nu} \right) ds(x) = 0.
			\end{equation}
			Define the function $w$ by
			\begin{equation}
				w(x)=\sum_{n\in\mathbb{Z}^{2}}v_{n}e^{i(\hat{\alpha}_{n}\cdot\tilde{x}+\hat{\beta}_{n}(2b-x_{3}))},\quad 0<x_{3}<b.
				\label{w(x)}
			\end{equation}
			Together with (\ref{G(x,y)}) and (\ref{w(x)}), we can compute
			\begin{align*}
				& \int_{\Gamma_{b}} \left( G(x,y)\frac{\partial w}{\partial\nu} - w\frac{\partial G(x,y)}{\partial\nu} \right) ds(x)\\
				&=\frac{i}{8\pi^{2}}\sum\limits_{n,m\in\mathbb{Z}^{2}}\frac{1}{\beta_{n}}v_{m}e^{-i\alpha_{n}\cdot\tilde{y}}e^{i\beta_{n}(y_{3}-b)}e^{i\hat{\beta}_{m}b}\cdot(-i\hat{\beta}_{m}+i
				\beta_{n})\int_{0}^{2\pi}\int_{0}^{2\pi}e^{i(\alpha_{n}+\hat{\alpha}_{m})\cdot\tilde{x}}d\tilde{x}\\
				&=\frac{i}{2}\sum\limits_{n\in\mathbb{Z}^{2}}\frac{1}{\beta_{n}}v_{-n}e^{-i\alpha_{n}\cdot\tilde{y}}e^{i\beta_{n}(y_{3}-b)}e^{i\hat{\beta}_{-n}b}\cdot(-i\hat{\beta}_{-n}+i
				\beta_{n})\\
				&=0.
			\end{align*}
			Here we have used the fact that $\hat{\beta}_{-n}=\beta_{n}$ and $\hat{\alpha}_{-n}+\alpha_{n}=0$. Since $v_{+}=w$ on $\Gamma_{b}$, we further deduce that
			\begin{equation}
				\int_{\Gamma_{b}} \left( G(x,y)\frac{\partial w}{\partial\nu} - w\frac{\partial G(x,y)}{\partial\nu} \right) ds(x)
				= \int_{\Gamma_{b}} \left( G(x,y)\frac{\partial w}{\partial\nu} - v_{+}\frac{\partial G(x,y)}{\partial\nu} \right) ds(x) = 0.
				\label{3.6}
			\end{equation}
			By (\ref{3.3})-(\ref{3.6}), we obtain
			\begin{equation}
				\begin{aligned}
					&-\int_{\Gamma_{b}} R_{\eta}u(x,y)\overline{\varphi}ds(x)\\
					&= \int_{\Gamma_{b}} \left( u(x,y)\frac{\partial v_{+}}{\partial\nu} - v_{+}\frac{\partial u(x,y)}{\partial\nu} \right) ds(x) \\
					&= \int_{\Gamma_{b}} \left( G(x,y)\frac{\partial v_{+}}{\partial\nu} - v_{+}\frac{\partial G(x,y)}{\partial\nu} \right) ds(x)+ \int_{\Gamma_{b}} \left( u^{s}(x,y)\frac{\partial v_{+}}{\partial\nu} - v_{+}\frac{\partial u^{s}(x,y)}{\partial\nu} \right) ds(x) \\
					&= \int_{\Gamma_{b}} G(x,y) \left( \frac{\partial v_{+}}{\partial\nu} - \frac{\partial w}{\partial\nu} \right) ds(x).
				\end{aligned}
				\label{3.7}
			\end{equation}
			Define the single-layer potential
			\[
			U(y)=\int_{\Gamma_{b}}G(x,y)h(x)ds(x),\quad y\in((0,2\pi)^{2}\times\mathbb{R})\backslash\Gamma_{b},
			\]
			where $h(x)=\frac{\partial v_{+}}{\partial\nu}-\frac{\partial w}{\partial\nu}\in H_{-\alpha}^{-\frac{1}{2}}(\Gamma_{b})$.
			From the assumption that $\int_{\Gamma_{b}}R_{\eta}u(x,y)\overline{\varphi}ds(x)=0$ and (\ref{3.7}), it follows that $U(y)=0$ for $y\in\Gamma_{a}$. Consequently, $U(y)$ satisfies the exterior Dirichlet problem:
			\begin{equation}
				\begin{cases}
					\Delta U + k^{2}U = 0 & \text{in}\quad (0,2\pi)^{2}\times(a,\infty), \\
					U=0 & \text{on}\quad \Gamma_{a},\\
					U(y)=\sum\limits_{n\in\mathbb{Z}^{2}}A^+_{n}e^{i(\hat{\alpha}_{n}\cdot\tilde{y}+\hat{\beta}_{n}(y_{3}-b))} & \text{for}\quad y_{3}>a.
				\end{cases}
				\label{2.7}
			\end{equation}
			By the uniqueness of solutions to the Dirichlet problem (\ref{2.7}) and the unique continuation principle, we deduce that $U(y)=0$ in the region $y_{3}>b$. This implies that $U_{+}(y)=\frac{\partial U_{+}(y)}{\partial\nu}=0$ on $\Gamma_{b}$. From the jump relation of the single-layer potential, we obtain $U_{-}(y)=U_{+}(y)=0$ on $\Gamma_{b}$. Hence, $U(y)$ satisfies the Dirichlet problem:
			\begin{equation}
				\begin{cases}
					\Delta U + k^{2}U = 0 & \text{in}\quad (0,2\pi)^{2}\times(-\infty,b), \\
					U=0 & \text{on}\quad \Gamma_{b},\\
					U(y)=\sum\limits_{n\in\mathbb{Z}^{2}}A^-_{n}e^{i(\hat{\alpha}_{n}\cdot\tilde{y}-\hat{\beta}_{n}(y_{3}-b))}& \text{for}\quad y_{3}<b.
				\end{cases}
				\label{2.8}
			\end{equation}
			Due to the uniqueness of downward solution to the problem (\ref{2.8}), we have $U(y)=0$ in $y_{3}<b$, thus $\frac{\partial U_{-}(y)}{\partial\nu}=0$ on $\Gamma_{b}$. Moreover, the jump relation gives
			\[
			h(x)=\frac{\partial U_{-}}{\partial\nu}|_{\Gamma_{b}}-\frac{\partial U_{+}}{\partial\nu}|_{\Gamma_{b}}=0,
			\]
			which implies $\frac{\partial v_{+}}{\partial\nu}-\frac{\partial w}{\partial\nu}=0$ on $\Gamma_{b}$. Then, from the definition of the functions $v$ and $w$, we obtain the following transmission problem
			\begin{equation}
				\begin{cases}
					\Delta v+k^{2}v=0 & \text{in}\quad \Omega_{b}^{c},\\
					\Delta w+k^{2}w=0 & \text{in}\quad \Omega_{b},\\
					v_{+}-w=0 & \text{on}\quad \Gamma_{b},\\
					\frac{\partial v_{+}}{\partial\nu}-\frac{\partial w}{\partial\nu}=0 & \text{on}\quad \Gamma_{b},
				\end{cases}
				\label{transmission problem}
			\end{equation}
			with the Rayleigh radiation condition (\ref{rayleigh expansion2}). By the uniqueness of the solution to the transmission problem (\ref{transmission problem}), we conclude that $v=0$ in $\Omega_{b}^{c}$ and $w=0$ in $\Omega_{b}$. Thus $v_{+}=\frac{\partial v_{+}}{\partial\nu}=0$ on $\Gamma_{b}$. Combining with the boundary condition of the problem (\ref{transmission problem1}), we get $R_{\eta}v_{-}=0$ on $\Gamma_{b}$. Furthermore, $v$ satisfies the Helmholtz equation $\Delta v+k^{2}qv=0$ in $\Omega_{b}$ with the Dirichlet boundary condition $v=0$ on $\Gamma_{0}$. By Lemma \ref{uniquessrobin}, we have $v=0$ in $\Omega_{b}$, which gives $v_{-}|_{\Gamma_{b}}=0$. Finally, we get $\overline{\varphi}=v_{+}|_{\Gamma_{b}}-v_{-}|_{\Gamma_{b}}=0$, which completes the proof of Lemma \ref{lemma2.1}.
		\end{proof}
		\label{lemma2.1}
	\end{lemma}
	
	With the aforementioned density result, we introduce an interior  RtD operator on $\Gamma_b$. Subsequently, we use coincidence of the near-field data of $u_{1}(x,y)=u_2(x,y)$ on $\Gamma_{a}$ to demonstrate that the Cauchy data of $u_{j}(x,y)$ must be identical on $\Gamma_{b}$. Here $u_{j} (x,y)(j=1,2)$ denotes the total field corresponding to the incident field $u^{in}(x,y)$ and the refractive index $q_{j} (j=1,2)$. Then, it follows that the interior RtD operators corresponding to the two refractive indices are identical.
	
	For any $f\in H_{\alpha}^{-\frac{1}{2}}(\Gamma_{b})$, the interior RtD operator is defined as follows:
	\[
	T:H_{\alpha}^{-\frac{1}{2}}(\Gamma_{b})\longrightarrow H_{\alpha}^{\frac{1}{2}}(\Gamma_{b})\quad\quad T(f):=  w|_{\Gamma_{b}},
	\]
	where $w\in H_{\alpha}^{1}(\Omega_{b})$ solves the quasi-periodic boundary value problem:
	\begin{equation}
		\begin{cases}
			\Delta w+k^{2}qw=0 & \text{in} \quad \Omega_{b},\\
			w=0 & \text{on} \quad \Gamma_{0},\\
			\partial_{\nu}w-i\eta w=f & \text{on} \quad \Gamma_{b}.\\
		\end{cases}
		\label{dtn2}
	\end{equation}
By Lemma \ref{uniquessrobin}, $T$ is a continuous operator from $H_{\alpha}^{-\frac{1}{2}}(\Gamma_{b})$ to $H_{\alpha}^{\frac{1}{2}}(\Gamma_{b})$. In the following discussion, we denote by $T_{j}$ the interior RtD operator corresponding to problem (\ref{dtn2}) with the refractive index $q$ replaced by $q_{j}~(j=1,2)$.
	
	\begin{lemma}
		For any $y\in\Gamma_{a}$, let $u_{j}(\cdot,y)~(j=1,2)$ denote the $\alpha$-quasiperiodic total field of the problem (\ref{1.1})-(\ref{1.4}) corresponding to the refractive index $q_{j}$ and the incident point source $G(x,y)$. Assume that the near-field data on the measurement plane $\Gamma_{a}$ coincide in the sense that
		\begin{equation}
			u_{1}(x,y)=u_{2}(x,y)\quad\text{for all} \quad x,y\in \Gamma_{a},x\ne y,
			\label{u1=u2}
		\end{equation}
		then
		\[
		T_1(f)=T_2(f)\quad\text{for all}\quad f\in H_{\alpha}^{-\frac{1}{2}}(\Gamma_{b}).
		\]
		\label{lemma2.3}
	\end{lemma}
	\begin{proof}
		Let $\tilde{u}(x,y)=u_{1}(x,y)-u_{2}(x,y)$ for $y\in\Gamma_{a}$. We obtain $\tilde{u}(x,y)=0$ for all $x\in\Gamma_{a}\backslash\{y\}$ from (\ref{u1=u2}), and the continuous extension of $\tilde{u}(\cdot,y)$ implies that it vanishes at $x=y$; hence, $\tilde{u}(x,y)=0$ on $\Gamma_{a}$. Since $q_{1}=q_{2}=1$ in the half-space $\{x_{3}>a\}$, the function $\tilde{u}$ satisfies the following boundary value problem:
		\begin{equation*}
			\begin{cases}
				\Delta\tilde{u}+k^{2}\tilde{u}=0 & \text{in} \quad (0,2\pi)^{2}\times(a,\infty),\\
				\tilde{u}=0 & \text{on} \quad \Gamma_{a},\\
				\tilde{u}=\sum\limits_{n\in\mathbb{Z}^{2}}\hat{u}_{n}(y)e^{i\alpha_{n}\cdot\tilde{x}}e^{i\beta_{n}(x_{3}-a)},& x_{3}\ge a.
			\end{cases}
		\end{equation*}
		By the uniqueness of the exterior Dirichlet problem and the unique continuation principle, we conclude that $\tilde{u}=0$ in the region $x_{3}\ge b$. This immediately yields:
		\begin{equation*}
			u_{1}^{+}(\cdot,y)|_{\Gamma_{b}}=u_{2}^{+}(\cdot,y)|_{\Gamma_{b}},\quad\partial_{\nu}u_{1}^{+}(\cdot,y)|_{\Gamma_{b}}=\partial_{\nu}u_{2}^{+}(\cdot,y)|_{\Gamma_{b}}\quad\text{for all}\quad y\in\Gamma_{a}.
			\label{u1(x,y)=u2(x,y)}
		\end{equation*}
		In addition, combining with the transmission conditions on $\Gamma_{b}$,
		\[
		u_{j}^{+}(\cdot,y)|_{\Gamma_{b}}=u_{j}^{-}(\cdot,y)|_{\Gamma_{b}},\quad\partial_{\nu}u_{j}^{+}(\cdot,y)|_{\Gamma_{b}}=\partial_{\nu}u_{j}^{-}(\cdot,y)|_{\Gamma_{b}}~\ (j=1,2),
		\]
		we have
		\begin{equation*}
			u_{1}(x,y)=u_{2}(x,y),\quad\partial_{\nu}u_{1}(x,y)=\partial_{\nu}u_{2}(x,y)\quad\text{for all $x\in\Gamma_{b},y\in\Gamma_{a}$}.
			\label{T1=T2}
		\end{equation*}
		Defining $f_{y}:=R_{\eta}u_{j}(\cdot,y)|_{\Gamma_{b}}\in H_{\alpha}^{-\frac{1}{2}}(\Gamma_{b})$ and using the definition of the interior RtD operators $T_{j}$, we further get
		\[
		T_{1}(f_{y})= u_{1}(\cdot,y)|_{\Gamma_{b}}= u_{2}(\cdot,y)|_{\Gamma_{b}}=T_{2}(f_{y})\quad\text{for every}\quad y\in\Gamma_{a}.
		\]
		By the linearity of $T_{j} (j=1,2)$, this identity can extend to all finite linear combinations of such data, that is
		\begin{equation}
			T_{1}(f)=T_{2}(f)\quad\text{for all } f\in \text{span}\{f_{y}:y\in\Gamma_{a}\}.
			\label{T1=T22}
		\end{equation}
		Write $\mathcal{G}:=\{R_{\eta}u_{j}(x,y)|_{\Gamma_{b}}\colon y\in\Gamma_{a}\}$. Then the span $\mathcal{G}$ is dense in $H_{\alpha}^{-\frac{1}{2}}(\Gamma_{b})$ according to Lemma \ref{lemma2.1}. Thus, for any $f\in H_{\alpha}^{-\frac{1}{2}}(\Gamma_{b})$, there exists a sequence $\{f_{n}\}_{n\in\mathbb{N}}\subset\text{span } \mathcal{G}$ such that
		\[
		\left \| f_{n}-f \right \|_{H_{\alpha}^{-\frac{1}{2}}(\Gamma_{b})} \to0\quad\text{as}\quad n\to\infty.
		\]
		By (\ref{T1=T22}), we yield $T_{1}(f_{n})=T_{2}(f_{n}) $ for each $n$. Using the continuity of $T_{j}$, we obtain
		\begin{equation*}
			\begin{aligned}
				\left \| T_{1}(f) -T_{2}(f)\right \|_{H_{\alpha}^{\frac{1}{2}}(\Gamma_{b})}&=\left \| T_{1}(f-f_{n}) -T_{2}(f-f_{n})\right \|_{H_{\alpha}^{\frac{1}{2}}(\Gamma_{b})}\\
				&\le\left \|  T_{1}(f-f_{n})\right \|_{H_{\alpha}^{\frac{1}{2}}(\Gamma_{b})}+\left \| T_{2}(f_{n}-f) \right \|_{H_{\alpha}^{\frac{1}{2}}(\Gamma_{b})}\\
				&\le C_{1}\left \| f-f_{n} \right \| _{H_{\alpha}^{-\frac{1}{2}}(\Gamma_{b})}+C_{2}\left \| f_{n}-f \right \| _{H_{\alpha}^{-\frac{1}{2}}(\Gamma_{b})}\\
				&\to0,\quad \text{as }n\to\infty,
			\end{aligned}
		\end{equation*}
		where constants $C_{1},C_{2}>0$ are independent of $n$. Hence,
		\[
		T_{1}(f)=T_{2}(f)\quad\text{for all}\quad f\in H_{\alpha}^{-\frac{1}{2}}(\Gamma_{b}),
		\]
		which proves $T_{1}=T_{2}$ as the operators from $H_{\alpha}^{-\frac{1}{2}}(\Gamma_{b})$ into $H_{\alpha}^{\frac{1}{2}}(\Gamma_{b})$.
	\end{proof}
	
	We can now use Lemma \ref{lemma2.3} and Green's formula to derive an orthogonality relation for the refractive indices.
	
	\begin{lemma}
		If $u_{1}(x,y)=u_{2}(x,y)$ for all $x,y\in \Gamma_{a},x\ne y$, then we obtain the following orthogonality relation:
		\begin{equation}
			\int_{\Omega_{b}} \phi_{1}(x) \overline{\phi_{2}(x)} \bigl[q_{1}(x) - q_{2}(x)\bigr] \,dx = 0,   \label{orthogonality relation}
		\end{equation}
		where $\phi_{1},\phi_{2}\in H_{\alpha,0}^{1}(\Omega_{b})$ are $\alpha$-quasiperiodic solutions of $\Delta\phi_{j}+k^{2}q_{j}\phi_{j}=0$ in $\Omega_{b}$ satisfying $\phi_{j}=0$ on $\Gamma_{0}$ ($j=1,2$).
		\begin{proof}
			Let $f=R_{\eta}\phi_{1}(x)|_{\Gamma_{b}}\in H_{\alpha}^{-\frac{1}{2}}(\Gamma_{b})$. Assume that $w_{1}\in H_{\alpha}^{1}(\Omega_{b})$ is the solution of the following boundary value problem:
			\begin{equation*}
				\begin{cases}
					\Delta w_{1}+k^{2}q_{2}w_{1}=0 &\text{in} \quad\Omega_{b},\\
					R_{\eta}w_{1}=f &\text{on}\quad\Gamma_{b},\\
					w_{1}=0 &\text{on}\quad\Gamma_{0}.
				\end{cases}
			\end{equation*}
			Define $\tilde{\phi}:=\phi_{1}-w_{1}$. Under the assumption that $u_{1}(x,y)=u_{2}(x,y)$, and combining with $T_{1}=T_{2}$ from Lemma \ref{lemma2.3}, we get $T_{1}f=\phi_{1}|_{\Gamma_{b}}=w_{1}|_{\Gamma_{b}}=T_{2}f$. Hence $\frac{\partial\phi_{1}}{\partial\nu}=f+i\eta\phi_{1}=f+i\eta w_{1}=\frac{\partial w_{1}}{\partial\nu}$, which implies that  $\frac{\partial \tilde{\phi}}{\partial\nu}=0$ on $\Gamma_{b}$. Consequently, $\tilde{\phi}$ satisfies the mixed boundary value problem:
			\begin{equation}
				\begin{cases}
					\Delta \tilde{\phi}+k^{2}q_{2}\tilde{\phi}=k^{2}(q_{2}-q_{1})\phi_{1} &\text{in}\quad\Omega_{b},\\
					\tilde{\phi}=0 &\text{on}\quad\Gamma_{0}\cup\Gamma_{b},\\
					\frac{\partial \tilde{\phi}}{\partial\nu}=0 &\text{on}\quad\Gamma_{b}.
				\end{cases}
				\label{q}
			\end{equation}
			Additionally, $\phi_{2}$ satisfies the Helmholtz equation $\Delta \phi_{2}+k^{2}q_{2}\phi_{2}=0$ in $\Omega_{b}$ (recall that $q_{2}$ is real valued) with the homogeneous Dirichlet boundary condition $\phi_{2}=0$ on $\Gamma_{0}$. To proceed, we multiply both sides of (\ref{q}) by $\overline{\phi_{2}}$ and integrate over $\Omega_{b}$. Applying Green's first identity and using the boundary conditions, we obtain:
			\begin{equation*}
				\begin{aligned}
					\int_{\Omega_{b}}k^{2}(q_{2}-q_{1})\phi_{1}\overline{\phi_{2}}dx&=\int_{\Omega_{b}}(\Delta \tilde{\phi}+k^{2}q_{2}\tilde{\phi})\overline{\phi_{2}}dx\\
					&=\int_{\Gamma_{0}\cup\Gamma_{b}}\frac{\partial \tilde{\phi}}{\partial\nu}\overline{\phi_{2}}ds-\int_{\Omega_{b}}\left[\nabla \tilde{\phi}\nabla\overline{\phi_{2}}-k^{2}q_{2}\tilde{\phi}\overline{\phi_{2}}\right]dx\\
					&=\int_{\Gamma_{0}\cup\Gamma_{b}}\left[\frac{\partial \tilde{\phi}}{\partial\nu}\overline{\phi_{2}}-\frac{\partial\overline{\phi_{2}}}{\partial\nu}\tilde{\phi}\right]ds\\
					&=0.
				\end{aligned}
			\end{equation*}
		By the quasi-periodic boundary condition, the contributions from the vertical sides of $\partial\Omega_{b}$ cancel out. Moreover, all remaining boundary terms vanish: on $\Gamma_{0}$ both $\tilde{\phi}$ and $\phi_{2}$ vanish, while on $\Gamma_{b}$ we have $\tilde{\phi}=0$ and $\partial_{\nu}\tilde{\phi}=0$ by (\ref{q}). 
		This completes the proof.
		\end{proof}
	\end{lemma}

	\section{High-Energy Plane-Wave Approximation via Spectral Projections}\label{Sec4}
	
	To establish the uniqueness result for the inverse problem, we need to construct a special solution of the original problem (\ref{1.1})-(\ref{1.4}) and substitute them into the orthogonality relation (\ref{orthogonality relation}) to recover the Fourier coefficients of the difference of the refractive indices. In two dimensions \cite{Kirsch A1995}, this construction relies on the high-energy spectral properties of an associated one-dimensional periodic Schr\"odinger operator. We are going to extend this approach to three dimensions for bi-periodic refractive indices. For this purpose,  it is necessary to introduce some basic knowledge about the Schr\"odinger operator with bi-periodic potential in $\mathbb{R}^{2}$ (\cite{V2024}).
	
	\subsection{Preliminary knowledge}
	
	Let $\mathcal{L}$ be a physical lattice in $\mathbb{R}^{2}$ with reduced basis $\omega_{1}=(2\pi,0), \omega_{2}=(0,2\pi)$, that is, $\mathcal{L}$$=\{\omega=n_{1}\omega_{1}+n_{2}\omega_{2}=(2\pi n_{1},2\pi n_{2})\in2\pi\mathbb{Z}^{2} \}$. We denote the dual lattice of $\mathcal{L}$ by $\mathcal{L}^{*}$, with basis $\mathfrak{e} _{j}$ biorthogonal to $\omega_{i}$ i.e. $\omega_{i}\cdot\mathfrak{e} _{j}=2\pi\delta_{ij}(i,j=1,2)$. Hence $\mathfrak{e} _{1}=(1,0),\mathfrak{e} _{2}=(0,1)$, and the reciprocal lattice is defined as $\mathcal{L}^{*}$$:=\{ n\in\mathbb{R}^{2}:n=n_{1}\mathfrak{e} _{1}+n_{2}\mathfrak{e} _{2}=(n_{1},n_{2})\in\mathbb{Z}^{2} \}$. The fundamental domain of the lattice $\mathcal{L}$ is taken as $F=\{ \tilde{x}\in \mathbb{R}^{2}:\tilde{x}=a_{1}\omega _{1}+a_{2}\omega_{2},0\le a_{i}<1 \}=[0,2\pi)^{2}$. 
	In addition, we usually consider the first Brillouin zone as $F^{*}=[-\frac{1}{2},\frac{1}{2}]^{2}$.
	
	We need to consider the following bi-periodic stationary Schr\"odinger equation in $\mathbb{R}^{2}$:
	\begin{equation}
		L(\tilde{q})\psi=[-\Delta_{\tilde{x}}+\tilde{q}(\tilde{x})]\psi(\tilde{x})=\Lambda \psi(\tilde{x}),\quad\tilde{x}\in\mathbb{R}^{2}.
		\label{L(q)}
	\end{equation}
	with $\tilde{q}(\tilde{x})=-k^{2}q(\tilde{x})$. The potential function $\tilde{q}(\tilde{x})\in L_{per}^{2}(F)$  is assumed to be real-valued and $2\pi$-periodic in each variable. Its Fourier expansion is written as:
	\begin{equation*}
		\tilde{q}(\tilde{x})=\sum\limits _{n\in\mathcal{L}^{*}}\tilde{q}_{n}e^{in\cdot\tilde{x}}=\sum\limits _{n\in \mathbb{Z}^2}\tilde{q}_{n}e^{in\cdot\tilde{x}},		
		\qquad \tilde{q}_{n}=\frac{1}{(2\pi)^{2}}\int_{F}\tilde{q}(\tilde{x})e^{-in\cdot\tilde{x}}d\tilde{x},
		\label{q1}
	\end{equation*}
	with $\tilde{q}_{-n}=\overline{\tilde{q}_{n}}$. 
For each quasi-momentum $t=(t_{1},t_{2})\in F^{*}$, we consider the eigenvalue problem (\ref{L(q)}) with the quasi-periodic boundary condition:
	\[
	\psi(\tilde{x}+\omega)=e^{it\cdot\omega}\psi(\tilde{x}),\quad\tilde{x}\in(0,2\pi)^{2},\omega\in\mathcal{L}.
	\]
	From Bloch's theorem, the eigenfunction $\psi$ of $L(\tilde{q})$ can be written in the Bloch form $\psi(\tilde{x})=e^{it\cdot\tilde{x} }p(\tilde{x})$, where the function $p(\tilde{x})$ is $2\pi$-biperiodic. Substituting this form into (\ref{L(q)}) yields the following eigenvalue problem for $p(\tilde{x})$ on $F$:
	\begin{equation}
		L_{t}(\tilde{q})p:=[-\Delta_{\tilde{x}}-2it\cdot\nabla_{\tilde{x}} +|t|^{2}+\tilde{q}(\tilde{x})]p(\tilde{x})=\Lambda p(\tilde{x}),
		\label{ltq}
	\end{equation}
	subject to the Floquet--Bloch boundary conditions:
	\begin{equation*}
		p(2\pi,x_{2})=p(0,x_{2}),~\ \frac{\partial p}{\partial x_{1}}(2\pi,x_{2})=\frac{\partial p}{\partial x_{1}}(0,x_{2});
	\end{equation*}
	\begin{equation*}
		p(x_{1},2\pi)=p(x_{1},0),~\ \frac{\partial p}{\partial x_{2}}(x_{1},2\pi)=\frac{\partial p}{\partial x_{2}}(x_{1},0).
	\end{equation*}
	The operator $L_t(\tilde{q})$ is defined in $L^2(F)$
	with domain $H_p^2(F)$, where $H_p^2(F)$ denotes the
	periodic Sobolev space of order two.
	Equivalently, in the $t$-quasiperiodic formulation,
	the corresponding operator $L(\tilde{q})=-\Delta+\tilde{q}$
	has domain
	\[
	H_t^2(F)
	=
	\left\{
	e^{it\cdot\tilde{x}}p(\tilde{x}):
	p\in H_p^2(F)
	\right\}.
	\]
	In both formulations, the function and its first derivatives
	satisfy the corresponding periodic or quasiperiodic matching
	conditions on opposite sides of the cell, understood in
	the trace sense.
	
	For each fixed $t\in F^{*}$, the operator $L_{t}(\tilde{q})$ is self-adjoint on $L^{2}(F)$ with compact resolvent. Hence, its spectrum consists of a discrete set of eigenvalues:
	\begin{equation*}
		\Lambda_{1}(t)\le\Lambda_{2}(t)\le\Lambda_{3}(t)\le\cdots\to+\infty,
	\end{equation*}
	where the eigenvalues are repeated according to multiplicity. Let $\{p_{m,t}\}_{m\in \mathbb{N}}$ be the eigenfunctions of $L_{t}(\tilde{q})$ corresponding to $\Lambda_{m}(t)$, and denote the associated Bloch function by $\psi_{m,t}=e^{it\cdot\tilde{x}}p_{m,t}(\tilde{x})$. Then $\{\psi_{m,t}\}_{m\in\mathbb N}$ is also a complete orthonormal
	family in $L^2(F)$ and
	\begin{equation*}
		[-\Delta_{\tilde{x}}+\tilde{q}(\tilde{x})]\psi_{m,t}
		=\Lambda_m(t)\psi_{m,t}.
	\end{equation*}
	Moreover, the spectrum of the periodic Schr\"odinger operator $L(\tilde{q})$ is given by
	\begin{equation*}
		\sigma(L(\tilde{q}))=\overline{\bigcup_{t\in F^{*}}\sigma (L_{t}(\tilde{q}))}
		=\overline{\bigcup_{m\in\mathbb{N}}\{ \Lambda_{m}(t),t\in F^{*}\}}.
	\end{equation*}
	
	In particular, for the unperturbed fiber operator, we have
	\[
	L_t(0)e^{in\cdot\tilde{x}}=|n+t|^2e^{in\cdot\tilde{x}},
	\qquad n\in\mathbb Z^2.
	\]
	Consequently, the high-energy regime corresponds to $|n+t|\to\infty$.
	In the following discussion, we retain all eigenfunctions in a suitable
	spectral interval around $|n+t|^2$. Below the inner product is denoted by
	$(f,g)_{L^2(F)}=\int_F f\overline g\,d\tilde{x}$ for any $f,g\in L^{2}(F)$.
	\subsection{High-energy asymptotic approximation}
	
	The following lemma provides the plane-wave approximation needed in
	the inverse problem. The sum includes the entire eigenspace whenever
	its eigenvalue belongs to the specified interval\cite{E. B. Davies1995}.
	
	\begin{lemma}[\textbf{Spectral projection}]\label{lem:r45-spectral-projection}
		Let $\tilde{q}$ be real valued, bounded, and $2\pi$-periodic in each
		variable, and let $t\in\mathbb R^2$. For a bounded closed interval
		$I\subset\mathbb R$, put
		\[
		V_{I,t}^{per}:=\operatorname{span}
		\{p_{m,t}:\Lambda_m(t)\in I\}.
		\]
		Then $V_{I,t}^{per}$ is finite dimensional, and the orthogonal projection
		$P_{I,t}=\mathbf{1}_{I}(L_{t}(\tilde{q}))$ onto this space is given by
		\begin{equation}\label{eq:r45-spectral-projection}
			P_{I,t}f
			=\sum_{\Lambda_m(t)\in I}
			(f,p_{m,t})_{L^2(F)}p_{m,t},
			\qquad f\in L^2(F).
		\end{equation}
		If $I=[\mu-\delta,\mu+\delta]$, where $\mu\in\mathbb R$ and
		$\delta>0$, then
		\begin{equation}\label{eq:r45-projection-residual}
			\|f-P_{I,t}f\|_{L^2(F)}
			\le\frac1\delta
			\|(L_{t}(\tilde{q})-\mu)f\|_{L^2(F)},
			\qquad f\in H^2_p(F).
		\end{equation}
	\end{lemma}
	
	\begin{proof}
		The spectral properties  give
		$\Lambda_m(t)\to+\infty$. Since $I$ is bounded, only finitely many
		indices satisfy $\Lambda_m(t)\in I$, so $V_{I,t}^{per}$ is finite
		dimensional. Put $a_m=(f,p_{m,t})_{L^2(F)}$ and 
		let $v:=\sum\limits_{\Lambda_{m}(t)\in I}a_{m}p_{m,t}$. Clearly, $v\in V_{I,t}^{per}$. For every $r$ such that $\Lambda_{r}(t)\in I$, orthonormality gives 
		\[
		\begin{aligned}
		(f-v,p_{r,t})_{L^{2}(F)}&=(f,p_{r,t})_{L^{2}(F)}-\sum\limits_{\Lambda_{m}(t)\in I}a_{m}(p_{m,t},p_{r,t})_{L^{2}(F)}\\
		&=a_{r}-\sum\limits_{\Lambda_{m}(t)\in I}a_{m}\delta_{mr}=0.
		\end{aligned}
		\]
		Hence $f-v$ is orthogonal to $V_{I,t}^{per}$.	Since this space is finite dimensional and therefore closed, $v$ is the unique orthogonal projection of $f$ onto $V_{I,t}^{per}$. 
		This proves \eqref{eq:r45-spectral-projection}.
	Since $\left\{p_{m,t}\right\}_{m\ge1}$ is a complete orthonormal family of $L^{2}(F)$, we have
	\[
	f=\sum_{m=1}^{\infty}a_mp_{m,t}=\sum\limits_{\Lambda_{m}(t)\in I}a_{m}p_{m,t}+\sum\limits_{\Lambda_{m}(t)\notin I}a_{m}p_{m,t}
	\qquad\text{in }L^2(F),
	\]
and Parseval's identity gives
	\begin{equation}\label{eq:r45-projection-tail}
					\|f-P_{I,t}f\|_{L^2(F)}^2
					=\sum_{\Lambda_m(t)\notin I}
					|(f,p_{m,t})_{L^2(F)}|^2.
				\end{equation}
		If the interval contains no eigenvalue, then $V_{I,t}^{per}=\{0\}$ and
		$P_{I,t}=0$.
Now let $f\in H^2_p(F)$ and $I=[\mu-\delta,\mu+\delta]$.
		By self-adjointness of $L_{t}(\tilde{q})$ yields,
	\begin{align*}
	\big((L_{t}(\tilde{q})-\mu)f,p_{m,t}\big)_{L^2(F)}
	&=\big(f,(L_{t}(\tilde{q})-\mu)p_{m,t}\big)_{L^2(F)}\\
	&=(f,\Lambda_m(t)p_{m,t}-\mu p_{m,t})\\
	&=(\Lambda_m(t)-\mu)a_m.
	\end{align*}
		Applying Parseval's identity to $(L_{t}(\tilde{q})-\mu)f$, we obtain
		\[
		\|(L_{t}(\tilde{q})-\mu)f\|_{L^2(F)}^2
		=\sum_{m=1}^{\infty}|\Lambda_m(t)-\mu|^2|a_m|^2.
		\]
		For $\Lambda_m(t)\notin I$, we have
	$|\Lambda_m(t)-\mu|>\delta>0$. It implies that $1\le\frac{|\Lambda_m(t)-\mu|^{2}}{\delta^{2}}$. 
Therefore,
\[
		\begin{aligned}
			\|f-P_{I,t}f\|_{L^2(F)}^2
			&=\sum_{\Lambda_m(t)\notin I}|a_m|^2\\
			&\le\frac1{\delta^2}
			\sum_{\Lambda_m(t)\notin I}|\Lambda_m(t)-\mu|^2|a_m|^2\\
			&\le\frac1{\delta^2}
			\sum_{m=1}^{\infty}|\Lambda_m(t)-\mu|^2|a_m|^2\\
			&=\frac1{\delta^2}
			\|(L_{t}(\tilde{q})-\mu)f\|_{L^2(F)}^2.
		\end{aligned}
		\]
		Taking square roots proves \eqref{eq:r45-projection-residual}.
	\end{proof}
	
	For $E_{n}(\tilde{x}):=e^{in\cdot\tilde{x}}$ and
	$\rho=|n+t|$, the exact identity
	\[
	(L_{t}(\tilde{q})-\rho^2)E_{n}=(|n|^{2}+2t\cdot n+|t|^{2}+\tilde{q}-|n+t|^{2})E_{n}=\tilde{q}E_{n}
	\]
	shows that the residual is bounded in $L^2(F)$ independently of $n$.
	Lemma~\ref{lem:r45-spectral-projection} therefore gives an
	approximation by the spectral projection onto an interval of width
	$2\rho$ centered at $\rho^2$. We record this consequence in the
	following lemma.
	
	\begin{lemma}\label{lem:r451-spectral-approximation}
		Let $\tilde{q}$ be a real-valued, bounded, $2\pi$-periodic potential.
		Fix $t\in\mathbb R^2$ and $n\in\mathbb Z^2$ such that
		$\rho:=|n+t|>1$. Set $E_{n}=e^{in\cdot\tilde{x}}$,
		$I_{n,t}=[\rho^2-\rho,\rho^2+\rho]$, and let
		$\Psi_{n,t}:=e^{it\cdot\tilde{x}}P_{I_{n,t},t}E_{n}$. Then
		\begin{equation}\label{eq:r451-window}
			\Psi_{n,t}:=
		\sum_{\Lambda_m(t)\in I_{n,t}}
		\big(e^{i(n+t)\cdot\tilde{x}},\psi_{m,t}\big)_{L^2(F)}\psi_{m,t}
		\end{equation}
		The function $\Psi_{n,t}$ is a finite combination of $t$-quasiperiodic
		eigenfunctions with positive eigenvalues. Moreover,
		\begin{equation}\label{eq:r451-plane-wave-approximation}
			\Psi_{n,t}(\tilde{x})
			=e^{i(n+t)\cdot\tilde{x}}+r_{n,t}(\tilde{x}),
			\qquad
			\|r_{n,t}\|_{L^2(F)}
			\le\frac{2\pi\|\tilde{q}\|_{L^\infty(F)}}{|n+t|}.
		\end{equation}
	\end{lemma}
	\begin{proof}
		Since $\psi_{m,t}=e^{it\cdot\tilde{x}}p_{m,t}$ and $(e^{in\cdot\tilde{x}},p_{m,t})=(e^{i(n+t)\cdot\tilde{x}},\psi_{m,t})$, by	\eqref{eq:r45-spectral-projection}, we get
	
		\[
		\begin{aligned}
			\Psi_{n,t}
		&=e^{it\cdot\tilde{x}}P_{I_{n,t},t}E_{n}=e^{it\cdot\tilde{x}}\sum\limits_{\Lambda_{m}(t)\in I}(e^{in\cdot\tilde{x}},p_{m,t})p_{m,t}\\
		&=
		\sum_{\Lambda_m(t)\in I_{n,t}}
		(e^{i(n+t)\cdot\tilde{x}},\psi_{m,t})_{L^2(F)}\psi_{m,t}.
		\end{aligned}
		\]
		It is the expression of the
		orthogonal projection in the chosen eigenfunction basis.
		Since $E_{n}\in H^2_p(F)$, we may apply
		\eqref{eq:r45-projection-residual} with
		$f=E_{n}$, $\mu=\rho^2$ and $\delta=\rho$. This gives
		\begin{equation}\label{eq:r451-tail-estimate}
			\begin{aligned}
				\|e^{i(n+t)\cdot\tilde{x}}-\Psi_{n,t}\|_{L^2(F)}&=\|E_{n}-P_{I_{n,t},t}E_{n}\|_{L^2(F)}\le\frac1\rho
				\|(L_{t}(\tilde{q})-\rho^2)E_{n}\|_{L^2(F)}\\
				&=\frac1\rho\|\tilde{q} E_{n}\|_{L^2(F)}
				\le\frac{2\pi\|\tilde{q}\|_{L^\infty(F)}}{\rho},
			\end{aligned}
		\end{equation}
		where $|E_{n}|=1$ and $|F|=(2\pi)^2$.
		Defining $r_{n,t}=\Psi_{n,t}-e^{i(n+t)\cdot\tilde{x}}$ proves
		\eqref{eq:r451-plane-wave-approximation}.
		The sum in \eqref{eq:r451-window} is finite by
		Lemma~\ref{lem:r45-spectral-projection}. For each retained index,
		\[
		\Lambda_m(t)\ge\rho^2-\rho=\rho(\rho-1)>0.
		\]
		Thus every eigenvalue appearing in this finite sum is positive.
		All its terms have quasiperiodicity $t$, as does $e^{it\cdot\tilde{x}}E_{n}$.
		Hence both $\Psi_{n,t}$ and $r_{n,t}$ have the same quasiperiodicity.
	\end{proof}
	
	\begin{remark}\label{rem:r451-fixed-alpha}
		The function $\Psi_{n,t}$ is a spectral projection of a plane wave
		and is generally a combination of eigenfunctions with different
		eigenvalues. Formula \eqref{eq:r451-plane-wave-approximation} is a
		leading-term asymptotic formula for this finite spectral combination.
		The remainder tends to zero in $L^2(F)$ as $|n+t|\to\infty$.
		In particular, for $\rho>\max\{1,\|\tilde{q}\|_{L^\infty(F)}\}$,
		\[
		\|\Psi_{n,t}\|_{L^2(F)}
		\ge2\pi-\frac{2\pi\|\tilde{q}\|_{L^\infty(F)}}{\rho}>0,
		\]
		so the projection is nonzero in the high-energy regime.
			The finite sum in \eqref{eq:r451-window} represents the spectral projection of
		the plane wave, with the approximation error controlled by \eqref{eq:r451-plane-wave-approximation}.
		
		In our inverse problem, the quasi-momentum is prescribed by $\alpha$.
		Writing $\alpha=n_0+t$ with $n_0\in\mathbb Z^2$ and $t\in F^*$
		only changes the labeling of Fourier modes. We may therefore take
		$t=\alpha$ throughout. In particular,
		Lemma~\ref{lem:r451-spectral-approximation} applies to every fixed real
		$\alpha$, including phases at which spectral multiplicities occur.
	\end{remark}

	\section{Proof of uniqueness}
	
	In this section, we state and prove the main uniqueness result of our inverse problem.
	
	\begin{theorem}
		Let $u_{j}(x,y)$ solve the problem (\ref{1.1})-(\ref{1.4}) with 
		$q=q_{j}\,  (j=1,2)$ corresponding to the incident wave $G(x,y)$, and suppose that each refractive index $q_{j}$ depends only on $\tilde{x}=(x_{1},x_{2})$. If $u_{1}(x,y)=u_{2}(x,y)$ for all $x,y\in \Gamma_{a} (x\ne y)$ with some $a>b$, then $q_{1}=q_{2}$ a.e. in $\mathbb{R}_{+}^{3}$.
		\label{uniqueness theorem}
	\end{theorem}
	\begin{proof}
		We seek $\alpha$-quasiperiodic solutions to the equation $\Delta \phi_{j}+k^2q_j \phi_{j}=0$ in the separated form:
		\begin{equation}
			\phi_{j}(\tilde{x},x_{3})=\psi^{(j)}(\tilde{x})h^{(j)}(x_{3}),\quad j=1,2,
			\label{phix()}
		\end{equation}
		where $\psi^{(j)}(\tilde{x})$ is $\alpha$-quasiperiodic. Substituting such $\phi_{j}(x)$ into the Helmholtz equation gives
		\[
		\Delta_{\tilde{x}}\psi^{(j)}(\tilde{x})h^{(j)}(x_{3})+h^{(j)^{''}}(x_{3})\psi^{(j)}(\tilde{x})+k^{2}q_{j}(\tilde{x})\psi^{(j)}(\tilde{x})h^{(j)}(x_{3})=0.
		\]
		To construct such solutions, it suffices to choose a nonzero function $\psi^{(j)}$ and a function $h^{(j)}$ satisfying
		\begin{equation*}
			\begin{cases}
				-\Delta_{\tilde{x}}\psi^{(j)}(\tilde{x})-k^{2}q_{j}(\tilde{x})\psi^{(j)}(\tilde{x})=\lambda^{(j)}\psi^{(j)}(\tilde{x}),\\
			h^{(j)\prime\prime}-\lambda^{(j)}h^{(j)}=0,\qquad h^{(j)}(0)=0,
			\end{cases}
		\end{equation*}
for some constant $\lambda^{(j)}\in\mathbb{R}$.
		
		For each potential function $q_{j}$, we construct a special kind of solution $\psi^{(j)}(\tilde{x})$ by considering the following eigenvalue problem on the fundamental cell $F=[0,2\pi)^{2}$:
		\begin{equation*}
			\begin{cases}
				-\Delta_{\tilde{x}}\psi^{(j)}(\tilde{x})-k^{2}q_{j}(\tilde{x})\psi^{(j)}(\tilde{x})=\lambda^{(j)}\psi^{(j)}(\tilde{x}),\quad\tilde{x}\in F,\\
				\psi^{(j)}(\tilde{x}) \text{ is $\alpha$-quasiperiodic}.
			\end{cases}
			\label{SturmLP}
		\end{equation*}
		By standard Floquet-Bloch theory (\cite{Karpeshina Y E.2006,Kuchment P A2012,ReedSimon1978}), we know that the spectrum of the fiber operator $L_{\alpha}(\tilde{q}_{j}):=(-i\nabla_{\tilde{x}}+\alpha)^{2}+\tilde{q}_{j}(\tilde{x})$ is discrete. For fixed quasi-momentum $\alpha$, denote its Bloch eigenpairs by
		\begin{equation*}
			\lambda_{m}^{(j)}(\alpha),\quad\psi_{m}^{(j)}(\tilde{x} ;\alpha),\quad m\in\mathbb{N}.
		\end{equation*} 
For every  eigenvalue $\lambda_{m}^{(j)}>0$, we consider the one-dimensional boundary value problem for $h_{m}^{(j)}(x_{3})$:
		\begin{equation}
			\begin{cases}
				h_{m}^{(j)^{''}}(x_{3})-\lambda_{m}^{(j)}h_{m}^{(j)}(x_{3})=0,\quad x_{3}\in(0,b),\\
				h_{m}^{(j)}(0)=0.
			\end{cases}
			\label{h(x3)}
		\end{equation}
	 The general solution of (\ref{h(x3)}) is
		\begin{equation}
			\begin{aligned}
				h_{m}^{(j)}(x_{3})=c_{m}(e^{\sqrt{\lambda _{m}^{(j)}}x_{3}}-e^{-\sqrt{\lambda _{m}^{(j)}}x_{3}})
				=2c_{m}\sinh(\sqrt{\lambda _{m}^{(j)}}x_{3}),\quad c_{m}\in\mathbb{C}\backslash\{0\}.
			\end{aligned}
			\label{hn}
		\end{equation}
		
	For any pair $\lambda_m^{(1)}>0$ and $\lambda_{m'}^{(2)}>0$, let $\phi_{1,m}(x)=\psi_{m}^{(1)}(\tilde{x})h_{m}^{(1)}(x_{3})$ be the solution constructed by (\ref{phix()}) corresponding to $q_{1}$ and $\phi_{2,m'}(x)=\psi_{m'}^{(2)}(\tilde{x})h_{m'}^{(2)}(x_{3})$ be the solution corresponding to $q_{2}$ (recall that $q_2$ is real valued). It follows from the orthogonality relation (\ref{orthogonality relation}) that
		\begin{equation}  0=\int_{\Omega_{b}}\phi_{1,m}(x)\overline{\phi_{2,m'}(x)}[q_{1}(x_{1},x_{2})-q_{2}(x_{1},x_{2})]dx=I_{m,m'}\cdot J_{m,m'},
			\label{IS cdot JS}
		\end{equation}
		where
		\begin{align}
			I_{m,m'}&:=\int_{0}^{2\pi}\int_{0}^{2\pi}(q_{1}-q_{2})\psi_{m}^{(1)}(x_{1},x_{2})\overline{\psi_{m'}^{(2)}(x_{1},x_{2})}dx_{1}dx_{2},
			\label{Ik} \\
		J_{m,m'}&:=\int_{0}^{b}h_{m}^{(1)}(x_{3})\overline{h_{m'}^{(2)}(x_{3})}dx_{3}.
			\label{Js}
		\end{align}
	Substituting the  expression \ref{hn} into (\ref{Js}), we obtain
		\begin{equation*}
			\begin{aligned}
				J_{m,m'}&=\int_{0}^{b}4c_{m}\overline{c_{m'}}\sinh(\sqrt{\lambda_{m}^{(1)}}x_{3})\sinh(\sqrt{\lambda_{m'}^{(2)}}x_{3})dx_{3}\\
				&=2c_{m}\overline{c_{m'}}\int_{0}^{b}\left[\cosh((\sqrt{\lambda_{m}^{(1)}}+\sqrt{\lambda_{m'}^{(2)}})x_{3})-\cosh(|\sqrt{\lambda_{m}^{(1)}}-\sqrt{\lambda_{m'}^{(2)}}|x_{3})\right]dx_{3}.
			\end{aligned}
		\end{equation*}
		Since the function $\cosh(x)$ is strictly increasing on $[0,\infty)$, we have $$\cosh((\sqrt{\lambda_{m}^{(1)}}+\sqrt{\lambda_{m'}^{(2)}})x_{3})>\cosh(|\sqrt{\lambda_{m}^{(1)}}-\sqrt{\lambda_{m'}^{(2)}}|x_{3})$$ for all $x_{3}\in(0,b)$. This implies that the integrand is strictly positive on the open interval $(0,b)$. Hence, $J_{m,m'}\ne0$ for $\lambda_m^{(1)}>0$ and $\lambda_{m'}^{(2)}>0$ , and (\ref{IS cdot JS}) gives $I_{m,m'}=0$. We next use the high-energy approximation in
		Lemma~\ref{lem:r451-spectral-approximation}.
		Put $\gamma_n=n+\alpha$ and $\rho_n=|\gamma_n|$.
		For $\rho_n>1$, define
		\[
		I_n=[\rho_n^2-\rho_n,\rho_n^2+\rho_n],\quad
			\mathcal M_n^{(j)}=\{m\in\mathbb N:\lambda_m^{(j)}\in I_n\},\quad
			c_{m,n}^{(j)}=\big(e^{i\gamma_n\cdot\tilde{x}},\psi_m^{(j)}\big)_{L^2(F)}.
		\]
		The spectral projection in \eqref{eq:r451-window}, applied to
		$\tilde{q}_j$ with $t=\alpha$, is therefore
		\begin{equation}\label{eq:r45-medium-projection}
			\Psi_n^{(j)}
			=\sum_{m\in\mathcal M_n^{(j)}}c_{m,n}^{(j)}\psi_m^{(j)},
			\qquad j=1,2.
		\end{equation}
		Here $n$ labels the plane wave and its spectral interval, whereas
		$m$ labels the eigenpairs of the corresponding operator.
		Lemma~\ref{lem:r451-spectral-approximation} gives
		\begin{equation}\label{eq:r451-medium-approximation}
			\begin{aligned}
				\Psi_n^{(j)}(\tilde{x})
				=e^{i\gamma_n\cdot\tilde{x}}+r_n^{(j)}(\tilde{x}),\quad
				\|r_n^{(j)}\|_{L^2(F)}
				\le\frac{2\pi k^2\|q_j\|_{L^\infty(F)}}{|\gamma_n|},
				\qquad j=1,2.
			\end{aligned}
		\end{equation}
		Each $\Psi_n^{(j)}$ is a finite combination of eigenfunctions
		with positive eigenvalues. Its corresponding  solution of the original problem in $\Omega_{b}$ is
		\[
		\Phi_{j,n}(\tilde{x},x_3)
		=\sum_{m\in\mathcal M_n^{(j)}}c_{m,n}^{(j)}
		\psi_m^{(j)}(\tilde{x})h_m^{(j)}(x_3).
		\]
		Having obtained $I_{m,m'}=0$ for every pair
		of positive eigenvalues, we can now take finite linear combinations.
		For any $\ell\in\mathbb Z^2$ with $\rho_n>1$ and $\rho_{n+\ell}>1$,
		we obtain
		\begin{equation}\label{eq:r451-combination-orthogonality}
				\int_F(q_1-q_2)\Psi_{n+\ell}^{(1)}
				\overline{\Psi_n^{(2)}}\,d\tilde{x}=
				\sum_{m\in\mathcal M_{n+\ell}^{(1)}}
				\sum_{m'\in\mathcal M_n^{(2)}}
				c_{m,n+\ell}^{(1)}\overline{c_{m',n}^{(2)}}I_{m,m'}=0.
		\end{equation}
 Fix an arbitrary index shift $\ell=(\ell_1,\ell_2)\in\mathbb Z^2$.
		Choose $n_s\in\mathbb{Z}^{2}$, $s\in\mathbb N$ such that $|n_{s}|\to\infty$ as $s\to\infty$. In particular, we may choose $n_{s}=(s,0)$. Set
		$\rho_s=|\gamma_{n_s}|$. Then
		\[
		\rho_s\longrightarrow\infty,
		\qquad
		\frac{|\gamma_{n_s+\ell}|}{\rho_s}\longrightarrow1
		\quad\text{as }s\to\infty.
		\]
		For all sufficiently large $s$, define
		\begin{equation}\label{eq:r45-Is}
			I_s:=\int_F(q_1-q_2)\Psi_{n_s+\ell}^{(1)}
			\overline{\Psi_{n_s}^{(2)}}\,d\\tilde{x}=0.
		\end{equation}
		By \eqref{eq:r451-medium-approximation}, we may write
		\begin{equation}\label{eq:r45-sequence-approximation}
			\begin{aligned}
				\Psi_{n_s+\ell}^{(1)}(\tilde{x})
				=e^{i\gamma_{n_s+\ell}\cdot\tilde{x}}+r_s^{(1)}(\tilde{x}),\quad
				\Psi_{n_s}^{(2)}(\tilde{x})
				=e^{i\gamma_{n_s}\cdot\tilde{x}}+r_s^{(2)}(\tilde{x}),
			\end{aligned}
		\end{equation}
		with $\left \|r_{s}^{(j)}  \right \|_{L^{2}(F)}\le C\rho_{s}^{-1} ,j=1,2$, where the constant $C>0$ is independent of $n\in\mathbb{Z}^{2}$. Substituting (\ref{eq:r45-sequence-approximation}) into (\ref{eq:r45-Is}), we can expand $I_{s}$ as
		\begin{equation*}
			\begin{aligned}
				I_{s}=&\int_{0}^{2\pi}\int_{0}^{2\pi}(q_{1}-q_{2})\left ( e^{i(n_{s}+\ell+\alpha)\cdot \tilde{x}} +r_{s}^{(1)}(\tilde{x}) \right )\cdot\left ( \overline{e^{i(n_{s}+\alpha)\cdot \tilde{x}} +r_{s}^{(2)}(\tilde{x}) } \right )  dx_{1}dx_{2},\\
				=&\int_{0}^{2\pi}\int_{0}^{2\pi}(q_{1}-q_{2})e^{i\ell\cdot\tilde{x}}dx_{1}dx_{2}+R_{s},
			\end{aligned}
		\end{equation*}
		where $R_{s}$ collects all terms containing at least one remainder $r_{s}^{(j)}$.
		Explicitly,
		$
		R_{s}=R_{s}^{(1)}+R_{s}^{(2)}+R_{s}^{(3)}
		$
		with
		\begin{equation*}
			\begin{aligned}
				R_{s}^{(1)}&=\int_{F}(q_{1}-q_{2})(\tilde{x})e^{i(n_{s}+\ell+\alpha)\cdot\tilde{x}}\overline{r_{s}^{(2)}(\tilde{x})}d\tilde{x},\\R_{s}^{(2)}&=\int_{F}(q_{1}-q_{2})(\tilde{x})r_{s}^{(1)}(\tilde{x})e^{-i(n_{s}+\alpha)\cdot\tilde{x}}d\tilde{x}, \\
				R_{s}^{(3)}&=\int_{F}(q_{1}-q_{2})(\tilde{x})r_{s}^{(1)}(\tilde{x})\overline{r_{s}^{(2)}(\tilde{x})}d\tilde{x}.
			\end{aligned}
		\end{equation*}
		Notice that the main term in $I_{s}$ is exactly the Fourier coefficient of $q_{1}-q_{2}$ at the frequency $\ell=(\ell_{1},\ell_{2})$:
		\[\int_{0}^{2\pi}\int_{0}^{2\pi}[q_{1}(x_{1},x_{2})-q_{2}(x_{1},x_{2})]e^{i\ell\cdot\tilde{x}}dx_{1}dx_{2}.
		\]
		The remainder $R_{s}$ can be estimated as follows. Since $q_{j}\in L^{\infty}(F)$, we have $|(q_{1}-q_{2})(\tilde{x})|\le\left \| (q_{1}-q_{2}) \right \|_{L^{\infty}}$ almost everywhere. For $R_{s}^{(1)}$, using H\"older's inequality yields
		\begin{equation}
			\begin{aligned}
				|R_{s}^{(1)}|&\le\int_{F}|(q_{1}-q_{2})(\tilde{x})||r_{s}^{(2)}(\tilde{x})|d\tilde{x}\le\left \| q_{1}-q_{2} \right \|_{L^{\infty}}\left \|  r_{s}^{(2)}(\tilde{x})\right \|_{L^{1}(F)}\\
				&\le C_{1}\left \|  r_{s}^{(2)}(\tilde{x})\right \|_{L^{2}(F)}(\int_{F}1^{2}d\tilde{x})^{\frac{1}{2}}\\
				&\le C_{2}\rho_{s}^{-1}.
			\end{aligned}
			\label{Rs1}
		\end{equation}
		Similarly,
		\begin{equation}
			|R_{s}^{(2)}|\le\int_{F}|(q_{1}-q_{2})(\tilde{x})||r_{s}^{(1)}(\tilde{x})|d\tilde{x}\le\left \| q_{1}-q_{2} \right \|_{L^{\infty}}\left \|  r_{s}^{(1)}(\tilde{x})\right \|_{L^{1}(F)}\le C_{3}\rho_{s}^{-1}.
			\label{Rs2}
		\end{equation}
		For $R_{s}^{(3)}$, we use the Cauchy--Schwarz inequality to estimate
		\begin{equation}
			\begin{aligned}
				|R_{s}^{(3)}|&\le\int_{F}|(q_{1}-q_{2})(\tilde{x})||r_{s}^{(1)}(\tilde{x})||r_{s}^{(2)}(\tilde{x})|d\tilde{x}\\
				&\le\left \| q_{1}-q_{2} \right \|_{L^{\infty}}\left \|  r_{s}^{(1)}(\tilde{x})\right \|_{L^{2}(F)}\left \|  r_{s}^{(2)}(\tilde{x})\right \|_{L^{2}(F)}\\
				&\le C_{4}\rho_{s}^{-2}\le C_{4}\rho_{s}^{-1}.
			\end{aligned}
			\label{Rs3}
		\end{equation}
		Combining the above estimates (\ref{Rs1})-(\ref{Rs3}), we obtain
		\[
		|R_{s}|\le|R_{s}^{(1)}|+|R_{s}^{(2)}|+|R_{s}^{(3)}|\le C\rho_{s}^{-1},\quad\text{as }s\to\infty,
		\]
		where $C=C_{2}+C_{3}+C_{4}$. From $\rho_{s}\to\infty$, we conclude that $|R_{s}|\le C\rho_{s}^{-1}\to0$ as $s\to\infty$.
		Since $I_{s}=0$ for all sufficiently large $s$, letting $s\to\infty$ in (\ref{eq:r45-Is}) gives
		\begin{equation*}
			\int_{F}[q_{1}(\tilde{x})-q_{2}(\tilde{x})]e^{i\ell\cdot\tilde{x}}d\tilde{x}=0\quad\text{for every $\ell\in\mathbb{Z}^{2}$}.
			\label{INT F}
		\end{equation*}
		Hence, all Fourier coefficients of the $2\pi$-periodic function $f:=q_{1}-q_{2}\in L^{2}(F)$ vanish. Since exponentials $\{e^{i\ell\cdot\tilde{x}}:\ell\in\mathbb{Z}^{2}\}$ form an orthogonal basis of the Hilbert space $L^{2}(F)$, we conclude $f=0$ a.e. in $F$. By periodicity,  this identity extends to the entire space $\mathbb{R}^{2}$, which completes the proof of $q_{1}=q_{2}$.
	\end{proof}
	
Finally, we record the following consequence, in which the data are generated by individual incident Rayleigh modes rather than by a point source.
	
	\begin{corollary}
		Let $u_{j,n}(j=1,2)$ be the total field associated with the refractive index $q_{j}=q_{j}(\tilde{x})$ corresponding to the downward incident Rayleigh mode $u_{n}^{i}(x)=e^{i\alpha_{n}\cdot\tilde{x}-i\beta_{n}x_{3}}$.
		If $u_{1,n}(x)|_{\Gamma_{b}}=u_{2,n}(x)|_{\Gamma_{b}}$ for all $n\in\mathbb{Z}^{2}$, then $q_{1}=q_{2}$.
	\end{corollary}
	\begin{proof}
		It suffices to prove the coincidence of the total fields on $\Gamma_{b}$ induced by each incident point source wave (\ref{G(x,y)}). Fix $a>b$. For any $y\in\Gamma_{a}$, the $\alpha$-quasiperiodic point source admits the Rayleigh expansion:
		\[
		G(x,y)=\sum\limits_{n\in\mathbb{Z}^{2}}c_{n}(y)u_{n}^{i}(x),\quad c_{n}(y)=\frac{i}{8\pi^{2}\beta_{n}}e^{-i\alpha_{n}\cdot\tilde{y}}e^{i\beta_{n}a}\quad x_{3}<a.
		\]
		Denote $f_{n}=(\frac{\partial u_{n}^{in}}{\partial x_{3}}-\mathbb{T}u_{n}^{in})|_{\Gamma_{b}}=-2i\beta_{n}e^{i\alpha_{n}\cdot\tilde{x}-i\beta_{n}b}$. Since $\beta_{n}=i\sqrt{|\alpha_{n}|^{2}-k^{2}}$ for large $|n|$, we have $$c_{n}(y)f_{n}(\tilde{x})=\frac{1}{4\pi^{2}}e^{i\alpha_{n}\cdot(\tilde{x}-\tilde{y})}e^{(b-a)\sqrt{|\alpha_{n}|^{2}-k^{2}}},$$ which decays exponentially as $|n|\to\infty$.  The input series therefore converges in $H_{\alpha}^{-\frac{1}{2}}(\Gamma_{b})$. By the linearity and continuous dependence of the direct scattering solution on the incident field, the total field generated by the point source satisfies the series expansion
		\[
		u_{j}(x,y)=\sum\limits_{n\in\mathbb{Z}^{2}}c_{n}(y)u_{j,n}(x).
		\]
		Combining with the assumption $u_{1,n}(x)|_{\Gamma_{b}}=u_{2,n}(x)|_{\Gamma_{b}}$, for any $y\in\Gamma_{a}$, we get
		\begin{equation*}
			u_{1}(\cdot,y)|_{\Gamma_{b}}-u_{2}(\cdot,y)|_{\Gamma_{b}}=\sum\limits_{n\in\mathbb{Z}^{2}}c_{n}(y)(u_{1,n}(x)|_{\Gamma_{b}}-u_{2,n}(x)|_{\Gamma_{b}})=0.
			\label{u1-u2}
		\end{equation*}
Since $q_{1}=q_{2}=1$ above $\Gamma_{b}$, the point-source total fields $u_{1}(\cdot,y)$ and $u_{2}(\cdot,y)$ satisfy the same Helmholtz equation in the half-space $\{x_{3}>b\}$ and have the same outgoing Rayleigh expansion. As they agree on $\Gamma_{b}$, the uniqueness of the exterior Dirichlet problem gives $u_{1}(\cdot,y)=u_{2}(\cdot,y)$ for $x_{3}>b$, and the unique continuation principle yields the coincidence of their Cauchy data on $\Gamma_{b}$. Hence $T_{1}=T_{2}$, and the orthogonality relation combined with the spectral argument of Theorem \ref{uniqueness theorem} implies $q_{1}=q_{2}$.
	\end{proof}
	
	\begin{remark}
		The present spectral framework suggests a possible extension to the time-harmonic Maxwell system with a bi-periodic refractive index. Such an extension would require a suitable vector valued counterpart of the periodic spectral analysis, taking into account the $\nabla\times\nabla\times$ structure and the associated boundary conditions. An alternative approach, both for the Helmholtz equation and the Maxwell system, is to construct $\alpha$-quasiperiodic complex geometrical optics solutions and then substitute them directly into an orthogonality identity to recover the Fourier coefficients. This would provide a complementary route to the spectral method developed here and will be investigated in future work.
	\end{remark}
	
	\vspace{1em}
	
	\textbf{Acknowledgement}
	\vspace{0.5em}
	
	The work of G.Hu is supported by the National Natural Science Foundation of China (No. 12425112), the Fundamental Research Funds for Central Universities in China (No. 63263073, 63261173) and a Key Program (No. 25JCZDJC00970) of Natural Science Foundation of Tianjin, China.

\end{document}